\UseRawInputEncoding
\documentclass[11pt]{article}
\usepackage{amsmath}
\numberwithin{equation}{section}
\usepackage{amsfonts,amssymb,amsmath}
\usepackage{mathrsfs}
\usepackage{color}
\newcommand{\beq}{\begin{equation}}
\newcommand{\enq}{\end{equation}}

\newtheorem{Theorem}{Theorem}[section]
\newtheorem{Lemma}[Theorem]{Lemma}
\newtheorem{Corollary}[Theorem]{Corollary}
\newtheorem{Definition}[Theorem]{Definition}
\newtheorem{Remark}[Theorem]{Remark}

\newcommand{\benu}{\begin{enumerate}}
\newcommand{\beqa}{\begin{eqnarray}}
\newcommand{\beqan}{\begin{eqnarray*}}
\newcommand{\eay}{\end{array}}
\newcommand{\edm}{\end{displaymath}}
\newcommand{\eenu}{\end{enumerate}}
\newcommand{\eeq}{\end{equation}}
\newcommand{\eeqa}{\end{eqnarray}}
\newcommand{\eeqan}{\end{eqnarray*}}

\newcommand{\br}{\begin{Remark}}
\newcommand{\er}{\end{Remark}}

\newcommand{\bqa}{\begin{eqnarray}}
\newcommand{\eqa}{\end{eqnarray}}
\newcommand{\bqw}{\begin{eqnarray*}}
\newcommand{\eqw}{\end{eqnarray*}}

\newcommand{\bea}{\begin{array}{cc}}
\newcommand{\ena}{\end{array}}

\allowdisplaybreaks[4]
\begin{document}
\begin{center}

{\large \bf Existence and upper semicontinuity of time-dependent attractors for the non-autonomous nonlocal diffusion equations}\\

 \vspace{0.20in}Bin Yang $^{1}$ $\ $  Yuming Qin$^{2,\ast}$\\
\end{center}
$^{1}$ College of Information Science and Technology, Donghua University, Shanghai, 201620, P. R. China.\\
$^{2}$ Department of  Mathematics, Institute for Nonlinear Science, Donghua University, Shanghai, 201620, P. R. China,
 \vspace{3mm}

\begin{abstract}
In this paper, under some appropriate assumptions, we prove the existence of the minimal time-dependent pullback $\mathcal D_{\sigma}^{\mathcal{H}_{t}}$-attractors ${\mathcal{A}}_{\mathcal D_{\sigma}^{\mathcal{H}_{t}}}$ for the non-autonomous nonlocal diffusion equations in time-dependent space $\mathcal{H}_{t}(\Omega)$. Next, in same phase space, using a priori estimate and energy methods we establish the existence of time-dependent pullback attractors $\left\{A_{\xi}(t)\right\}_{t \in \mathbb R}$ and the upper semicontinuity of
$\left\{A_{\xi}(t)\right\}_{t \in \mathbb R}$ and the global attractor $A$ of equation (\ref{1.1-2}) with $\xi = 0$, that is,
$$
\lim _{\xi \rightarrow 0^{+}} \operatorname{dist}_{\mathcal H_{t}}\left(A_{\xi}(t), A \right)=0.
$$
\end{abstract}
\hspace{10mm}{\bf 2020 MSC:} 35B40, 35B41, 35B65, 35K57.

\hspace{4mm}{\bf Keywords:} Non-autonomous nonlocal diffusion equations, global attractor, the minimal time-dependent pullback $  {\mathcal D_{\sigma}^{\mathcal{H}_{t}}}$-attractors, pullback attractors, upper semicontinuity.

\section{Introduction}
\setcounter{equation}{0}

\let\thefootnote\relax\footnote{*Corresponding author.}
\let\thefootnote\relax\footnote{E-mails: yuming$\_$qin@hotmail.com (Y. Qin), binyangdhu@163.com (B. Yang).}
\quad
In the past few decades, many scholars have devoted to obtaining the well-posedness of solutions by studying the attractors of partial differential equations (see \cite{chm2, cv, lm, r, yqlm.2}).
It is worth mentioning that diffusion equations are also a prevalent research direction since they are applied in many disciplines such as physics, chemistry and biology (see \cite{acr.2,ww.2,zb.2,zl.2}).

In this paper, we consider the following non-autonomous nonlocal diffusion equations
\begin{equation}
\left\{\begin{array}{ll}
u_{t}-\varepsilon(t) \Delta u_{t}-a(l(u)) \Delta u=f(u)+\xi h(x,t) & \text { in } \Omega \times(\tau, \infty), \\
u=0 & \text { on } \partial \Omega\times(\tau, \infty), \\
u(x, \tau)=u_{\tau}  &\,\, x \in \Omega,
\end{array}\right.\label{1.1-2}
\end{equation}
in time-dependent space $\mathcal H_{t}(\Omega)$, where $\Omega  \subset \mathbb{R}^{N}(N \ge 3)$ is a bounded domain with smooth boundary $\partial \Omega$, the initial time $\tau  \le t \in \mathbb R$, $\xi$ is a small positive parameter and the definition and properties of $\mathcal H_{t}(\Omega)$ can be seen in $\S 2$.

Assume the function $\varepsilon(t) \in C^{1}(\mathbb{R})$ is a decreasing bounded function with respect to the parameter $t$ satisfying
\begin{equation}
\lim _{t \rightarrow+\infty} \varepsilon(t)=0,
\label{1.2-2}
\end{equation}
and there exists a constant $L>0$ such that
\begin{equation}
\sup _{t \in \mathbb{R}}(|\varepsilon(t)|+|\varepsilon^{\prime}(t)|) \leq L.
\label{1.3-2}
\end{equation}

Besides, the function $a(l(u)) \in C\left(\mathbb{R}; \mathbb{R}^{+}\right)$ is the nonlocal diffusion term of equation $(\ref{1.1-2})$ and satisfying
\begin{equation}
0 <m \leqslant a(s) \leqslant M,  \quad \forall \, s \in \mathbb{R},
\label{1.4-2}
\end{equation}
where $m$ and $M$ are constants.
Furthermore, assume $l(u): L^{2}(\Omega) \to \mathbb R$ is a continuous linear functional acting on $u$ that satisfies for some $g \in L^{2}(\Omega)$,
\begin{equation}
l(u)=l_{g}(u)=\int_{\Omega} g(x) u(x) d x.
\label{Lg-2}
\end{equation}

In addition, suppose the nonlinear term $f \in C(\mathbb{R}, \mathbb{R})$ and satisfies
\begin{equation}
|f(u)-f(v)| \leq C\left(|u|^{p}+|v|^{p}+1\right)|u-v|,
\label{1.5-2}
 \end{equation}
and
\begin{equation}
\limsup _{|s| \rightarrow \infty} \frac{f(s)}{s}<\lambda_{1},
\label{1.6-2}
\end{equation}
where $p \le \frac{4}{{N - 2}}$, $u, v, s \in \mathbb R$ and $\lambda_{1}$ is the first eigenvalue of $-\Delta$ in $H_{0}^{1}(\Omega)$ with the homogeneous Dirichlet boundary conditions.

Throughout this paper, the inner product of $L^{2}(\Omega)$ is denoted by $(\cdot, \cdot)$, and the corresponding norm is written as ${\|\cdot\|}$. Then we also assume that the function $f$ admits the following decomposition
\begin{equation}
f=f_{0}+f_{1},
\label{1.7-2}
\end{equation}
with $f_{0}, f_{1} \in C(\mathbb{R}, \mathbb{R})$ satisfying
\begin{equation}
\left|f_{0}(u)\right| \leq C\left(|u|^{p+1}+|u|\right),
\label{1.8-2}
\end{equation}
\begin{equation}
({f_0}(u),u) \le (2a(l(u)) - {\varepsilon ^\prime }(t) - \frac{{{\lambda _1}}}{{{\lambda _1} + 1}}\varepsilon (t)){\left\| {\nabla u} \right\|^2} - \frac{{{\lambda _1}}}{{{\lambda _1} + 1}}{\left\| u \right\|^2}
\label{1.9-2}
\end{equation}
and
\begin{equation}
\left|f_{1}(u)\right| \leq C\left(|u|^{p+1}+1\right),
\label{1.10-2}
\end{equation}
for any $u \in \mathbb R$. Also assume the external force $h(x,t) \in L_{\mathrm{loc}}^{2}(\mathbb R ; L^{2}(\Omega))$.

Now we recall some related works of problem $(\ref{1.1-2})$. The literature review based on the classic general diffusion equation ${u_t} - \Delta u = f(u)$ and attractors in time-dependent space $\mathcal H_{t}(\Omega)$ can be seen in our paper \cite{qy.2}. The equation (\ref{1.1-2}) is a nonlocal diffusion equation and the function $a(\, \cdot \,)$ is a nonlocal diffusion term and many people considered analogous problems in the past decade. Caraballo, Herrera-Cobos and Mar\'{i}n-Rubio \cite{chm2} studied the existence of minimal pullback attractors of $u_{t} -a(l(u)) \Delta u=f(u)+h(t)$ in $L^{2}(\Omega)$ and $H^{1}(\Omega)$. Later on, Peng, Shang and Zheng \cite{psz} proved a similar result of this problem in $H_0^1(\Omega )$. In addition, Caraballo, Herrera-Cobos and Mar\'{i}n-Rubio \cite{chm16} discussed the existence of pullback attractors of $u_{t}-(1-\varepsilon) a(l(u)) \Delta u=f(u)+\varepsilon h(t)$ with $\varepsilon  \in [0,1]$ and the upper semicontinuity of attractors in $L^{2}(\Omega)$, which means the family of pullback attractors corresponding global compact attractor associates with the autonomous nonlocal limit problem when $\varepsilon  \to 0$. Then they also proved the existence of the minimal pullback attractors for $p$-Laplacian reaction-diffusion equation $u_{t}-a(l(u)) \Delta_{p} u=f(u)+h(t)$ in $L^{2}(\Omega)$ and $L^{p}(\Omega)$ with $p \geq 2$ in \cite{chm3.2}. Besides, they established the existence of regular pullback attractors as well as their upper semicontinuous in $H^{1}(\Omega)$ of $u_{t}-g_{1}(\varepsilon) a(l(u)) \Delta u=\tilde{g}_{1}(\varepsilon) f(u)+g_{0}(\varepsilon) h(t)$ in \cite{chm4.2}, where $g_{1} \in C([0,1] ;(0, \infty)), \tilde{g}_{1} \in C([0,1] ;[0, \infty))$ and $g_{0} \in C([0,1])$. Additionally, in \cite{chm} they proved the existence of pullback $\mathcal D$-attractors in $L^{2}(\Omega)$ and $H^{1}(\Omega)$ .

Next we will introduce some results related to upper semicontinuity of diffusion equations. Guo and Wang \cite{gw.2} proved the global attractor $A$ of $u_{t}-\nu \Delta u+f(u)+\lambda_{0} u+g(x)=0$ is upper semicontinuity at 0 with respect to the global attractor $\{ {A_L}\} $, where $A$ and $\{ {A_L}\} $ are obtained when $\Omega  = \mathbb R$ and $\Omega  = [ - L,L]$, respectively. Later Carvalho,  Jos$\acute{e}$ and Robinson \cite{cjr1.2} took into account the continuity of pullback attractors for evolution processes. Furthermore, Wang and Qin \cite{wq.2} studied the upper semicontinuity of attractors of $u_{t}-\Delta u_{t}-\Delta u=f(u)+\varepsilon g(x, t)$ in $H_0^1(\Omega )$. Anh and Bao \cite{ab1.2} established the upper semicontinuity of pullback $\mathcal D$-attractors of $u_{t}-\varepsilon \Delta u_{t}-\Delta u+f(u)=g(t)$ in $H_0^1(\Omega )$. Besides, Anh and Bao \cite{ab2.2} and Wang \cite{w.2} obtained similar results of $u_{t}-\varepsilon \Delta u_{t}-\Delta u+f(x, u)+\lambda u=g(x, t)$ in $L^{2}(\mathbb{R}^{N})$ and $u_{t}-\varepsilon \triangle u_{t}-\triangle u=f(u)+g(x)+\varepsilon h(t)$ in $H_0^1(\Omega )$, respectively. Moreover, some authors also considered the upper semicontinuity between global attractors and uniform attractors (see \cite{ct2.2, ct1.2, wlq.2, wwq.2, xz.2}).

When the problem $(\ref{1.1-2})$ is compared with the general diffusion equations, it can be found that it adds the terms $a(\, \cdot \,)$, $\varepsilon(t)$, $f(u)$ and $\xi h(x, t)$, which undoubtedly increases challenges. Meanwhile, these terms make some conventional methods lose their effect. Therefore, we now introduce these difficulties and explain our methods, they should be creative fresh attempts.

$(1)$ The results of this paper are all obtained on the time-dependent space $\mathcal H_{t}(\Omega)$.
Thus since the time-dependent function $\varepsilon(t)$ exists in the norm of $\mathcal H_{t}(\Omega)$, when $u$ or $u_{t}$ is used as the test function of the equation $(\ref{1.1-2})_{1}$ in the energy estimation, the resulting equation cannot be directly estimated using the Gronwall inequality. To make it effective, we use the transformation
$$
\varepsilon(t) \frac{d}{d t}\|\nabla u\|^{2}=\frac{d}{d t}\left(\varepsilon(t)\|\nabla u\|^{2}\right)-\varepsilon^{\prime}(t)\|\nabla u\|^{2},
$$
which combined with $\varepsilon(t)$ is a decreasing function can prove the desired results.

$(2)$ In this paper, (\ref{1.2-2}) and (\ref{1.4-2}) are weaker than the conditions in Qin and Yang \cite{qy.2}. In fact, the weakening of these conditions forces us to be more refined in the calculations of the pullback absorbing sets. In addition, we use two methods to discuss the existence of the minimal pullback $\mathcal D_{\sigma}^{\mathcal{H}_{t}}$-attractors ${\mathcal{A}}_{\mathcal D_{\sigma}^{\mathcal{H}_{t}}}$ and pullback attractors $\mathcal{A}_{\xi}$. On the one hand, in $\S 4$ we prove the compactness of ${\mathcal{A}}_{\mathcal D_{\sigma}^{\mathcal{H}_{t}}}$ by introducing continuous functions. On the other hand, in $\S 5$ we use the theorem proved in Wang and Qin \cite{wq.2} to derive the desired results, which is also novel when compared with Qin and Yang \cite{qy.2}.

$(3)$ It is worth mentioning that the nonlinear term $f(u)$ and the external force term $\xi h(x, t)$ makes problem (\ref{1.1-2}) be studied in a more general functional framework. To obtain the dissipative properties of the process, we assume that $f$ satisfies $(\ref{1.7-2})-(\ref{1.10-2})$, which are weaker than the conditions in Wang and Qin \cite{wq.2}.

The structure of this paper is organized as follows. In $\S 2$, we introduce some useful abstract definitions, theorems and lemmas. Then in $\S 3$, we obtain the existence and uniqueness of the solutions of problem $(\ref{1.1-2})$ by the standard Faedo-Galerkin approximations and compactness argument. Furthermore, we shall derive the existence of the minimal pullback $\mathcal D_{\sigma}^{\mathcal{H}_{t}}$-attractors ${\mathcal{A}}_{\mathcal D_{\sigma}^{\mathcal{H}_{t}}}$, the pullback attractors $\mathcal{A}_{\xi}=\{{A}_{\xi}(t)\}_{t \in \mathbb R}$ and the upper semicontinuity of $\{{A}_{\xi}(t)\}_{t \in \mathbb R}$ and the global attractor $A$ of equation $(\ref{1.1-2})$ with $\xi=0$ in $\S 4$ and $\S 5$, respectively.

\section{Preliminaries}
$\ \ \ \ $ Before proving the main results, we first introduce some necessary abstract concepts in this section, such as basic definitions and properties of function spaces and attractors.

Let $\left\{X_{t}\right\}_{t \in \mathbb{R}}$ be a family of normed spaces with norm $\|\cdot\|_{X_{t}}$ and metric $d_{X_{t}}(\cdot, \cdot)$.
The closed $R$-ball with the origin as the center and $R$ as the radius in $\left\{X_{t}\right\}_{t \in \mathbb{R}}$ is denoted as $${\bar B_{{X_t}}}(0, R) = \left\{ {u \in {X_t}:\left\| u \right\|_{{X_t}}^2 \le R} \right\}.$$

The Hausdorff semidistance between two nonempty  sets $S_{1}, S_{2} \subset \left\{X_{t}\right\}_{t \in \mathbb{R}}$ is denoted by
$$
dist_{{X}_{t}}(S_{1},S_{2})=\sup _{x \in S_{1}} \inf _{y \in S_{2}}\|x-y\|_{{X}_{t}} \, .
$$

Next, we recall some notations about the family of the Hilbert spaces $D(A^{\frac{s}{2}})$ with $A=-\Delta$ and $s \in \mathbb{R}$. These spaces were widely used to study attractors (see \cite{ps.2, sy, wq.2}).
To simplify the notations, let $\mathcal {H}^{s}=D(A^{\frac{s}{2}})$ and its inner product and norm are denoted as ${( \cdot , \cdot )_{{{\cal H}^s}}} = ({A^{\frac{s}{2}}} \cdot ,{A^{\frac{s}{2}}} \cdot )$ and $\|\cdot\|_{\mathcal{H}^{s}}=\|A^{\frac{s}{2}}\cdot\|$, respectively.

\begin{Lemma} {\rm(\cite{ps.2})}\label{1m2.1-2} The properties of the space $\mathcal {H}^{s}=D(A^{\frac{s}{2}})$ are as follows:

$(i)$ Assume that $s_{1}>s_{2}$, then the embedding $D(A^{{\frac{s_{1}}{2}}}) \hookrightarrow D(A^{{\frac{s_{2}}{2}}})$ is compact.

$(ii)$ Assume that $s \in [0,\frac{n}{2})$, then the embedding $D(A^{{\frac{s}{2}}}) \hookrightarrow L^{\frac{2n}{n-2s}}(\Omega)$ is continuous.

$(iii)$ Assume that $s_{0}>s_{1}>s_{2}$, then for any $\delta > 0$, there exists constants $\epsilon$ and $C(\epsilon)=C_{\epsilon}(s_{0}, s_{1}, s_{2})$ such that
$$
\|A^{\frac{s_{1}}{2}}u\|\leq \epsilon\|A^{\frac{s_{0}}{2}}u\| + C(\epsilon)\|A^{\frac{s_{2}}{2}}u\|.
$$

$(iv)$ Assume that $s_{1}, s_{2} \in (0,1)$ and let $u \in \mathcal{H}^{s_{1}}(\Omega) \cap \mathcal{H}^{s_{2}}(\Omega)$, then for any constant $\theta \in (0,1)$, there exists a constant $C(\theta)>0$ such that
$$
\|u\|_{\mathcal{H}^{(1-\theta) s_{1} + \theta s_{2}}} \leq C(\theta)\|u\|_{\mathcal{H}^{s_{1}}}^{1-\theta}\, \|u\|_{\mathcal{H}^{s_{2}}}^{\theta}.
$$
\end{Lemma}
$\ \ \ \ $ In addition, for any $t \in \mathbb{R}$, the time-dependent space $\mathcal{H}_{t}(\Omega)$ is endowed with the norm
$$
\|u\|_{\mathcal{H}_{t}}^{2}=\|u\|_{2}^{2}+\varepsilon(t)\|\nabla u\|_{2}^{2}\,.
$$

In particular, assume the space ${\cal H}_t^{1+\alpha}(\Omega)$, more regular than $\mathcal{H}_{t}(\Omega)$, is endowed with the norm
$$
\|u\|_{\mathcal{H}_{t}^{1+\alpha}}^{2}=\| A^{\frac{\alpha }{2}}u\|_{2}^{2}+\varepsilon(t)\|{A^{\frac{{1{\rm{ + }}\alpha }}{2}}} u\|_{2}^{2}\,.
$$

\begin{Definition} {\rm(\cite{mwx,zxz})}
A process or a two-parameter semigroup on $\mathcal H_{t}(\Omega)$ is a family $\{U(t, \tau) \mid t, \tau \in \mathbb{R}, t \geqslant \tau\}$ of mapping $U(t, \tau): \mathcal H_{\tau} \rightarrow \mathcal H_{t}$ satisfies that $U(\tau, \tau)u=u$ for any $u \in \mathcal H_{\tau}$ and $U(t, s) U(s, \tau)=U(t, \tau)$ for all $t \geq s \geqslant \tau$.
\label{def2.1-2}
\end{Definition}

\begin{Definition} {\rm(\cite{mwx,zxz})}
The process ${\{ U(t,\tau )\} _{t \ge \tau }}$ on $\mathcal H_{t}(\Omega)$ is said to be continuous, if for any $t \ge \tau $ the mapping $U(t, \tau): \mathcal H_{\tau} \rightarrow \mathcal H_{t}$ is continuous.
\label{def2.2-2}
\end{Definition}

\begin{Definition} {\rm(\cite{mwx,zxz})}
The process ${\{ U(t,\tau )\} _{t \ge \tau }}$ on $\mathcal H_{t}(\Omega)$ is said to be closed, if for any sequence ${\{x_{n}\}\in \mathcal H_{t}(\Omega)}$ the equality $U(t,\tau)x=y$ can be concluded from $x_{n} \rightarrow x \in \mathcal H_{t}(\Omega)$ and $U(t, \tau) x_{n} \rightarrow y \in \mathcal H_{t}(\Omega)$.
\label{def2.3-2}
\end{Definition}

\begin{Remark}{\rm(\cite{bv.2})}\label{re2.1-2}
It is obvious that if a process is continuous, then it is closed.
\end{Remark}
\begin{Definition} {\rm(\cite{chm,psz})}
For any $\sigma>0$, let $\mathcal D$ be a nonempty class of all families of parameterized sets $\widehat{D}=\left\{D(t): t \in \mathbb{R}\right\} \subset \Gamma(\mathcal H_{t})$ such that
$$
\lim _{\tau \rightarrow-\infty}\left(e^{\sigma \tau} \sup _{u \in D(\tau)}\|u\|_{\mathcal H_{t}}^{2}\right)=0,
$$
where $\Gamma(\mathcal H_{t})$ denotes the family of all nonempty subsets of $\mathcal H_{t}(\Omega)$, then $\mathcal D$ will be called a tempered universe in $\Gamma(\mathcal H_{t})$.
\label{deftem-2}
\end{Definition}

\begin{Definition} {\rm(\cite{psz,zs2})}
A process ${\{ U(t,\tau )\} _{t \ge \tau }}$ is pullback $\mathcal D$-asymptotically compact on $\mathcal H_{t}(\Omega)$, if for any $t \in \mathbb{R}$, any ${\widehat D} \in {{\cal D}}$, any sequence ${\left\{ {{\tau _n}} \right\}_{n\in{\mathbb{N}^ + }}}\subset( - \infty ,t]$ and any sequence ${\left\{ {{x_n}} \right\}_{n \in {\mathbb{N}^ + }}} \subset {D({\tau}_{n})} \subset {{\cal H}_t}(\Omega )$, the sequence $\left.\left\{U (t, \tau\right) x_{n}\right\}_{n\in{\mathbb{N}^ + }}$ is relatively compact in $\mathcal H_{t}(\Omega)$ when ${\tau _n} \to  - \infty $.
\label{def2.4-2}
\end{Definition}

\begin{Definition} {\rm(\cite{psz,zs2})}
A family ${\widetilde D=\{\widetilde D(t): t \in \mathbb{R}\}}$ is pullback absorbing for the process ${\{ U(t,\tau )\} _{t \ge \tau }}$ on $\mathcal H_{t}(\Omega)$, if for any $t \in \mathbb{R}$ and any bounded subsets $B \subset \mathcal H_{t}(\Omega)$, there exists some constants $T(t, B)>0$ such that $U(t, t-\tau) B \subset \widetilde D(t)$ for any $\tau \geq T(t, B).$
\label{def2.5-2}
\end{Definition}

\begin{Definition} {\rm(\cite{psz,zs2})}
A family ${\widehat D_{0}=\left\{D_{0}(t): t \in \mathbb{R}\right\}} \in \Gamma(\mathcal H_{t})$ is pullback $\mathcal D$-absorbing for the process ${\{ U(t,\tau )\} _{t \ge \tau }}$ on $\mathcal H_{t}(\Omega)$, if for any $t \in \mathbb{R}$ and ${\widehat D} \in {{\cal D}}$, there exists a ${\tau _0} = {\tau _0}(t,{\widehat D}) < t$ such that
$
U(t, \tau) D(\tau) \subset D_{0}(t)
$
for any $\tau \leqslant \tau_{0}(t, \widehat{D})$.
\label{def2.6-2}
\end{Definition}

\begin{Definition} {\rm(\cite{cc})}
The set $A$ is called a global attractor for the  process ${\{ U(t,\tau )\} _{t \ge \tau }}$ on $\mathcal H_{t}(\Omega)$ if the following properties hold:

(i) $A$ is compact;

(ii) $A$ is invariant; and

(iii) $A$ attracts each bounded subset of $\mathcal H_{t}(\Omega)$.
\label{def2.7-2}
\end{Definition}

\begin{Definition} {\rm(\cite{wq.2})}
A family of compact sets $\mathcal{A}_{\xi}=\{{A}_{\xi}(t)\}_{t \in \mathbb R}$ is said to be a pullback attractor for the  process ${\{ U(t,\tau )\} _{t \ge \tau }}$ on $\mathcal H_{t}(\Omega)$ if the following properties hold:

(i) $\mathcal{A}_{\xi}$ is invariant, i.e., $U(t, \tau) {A}_{\xi}(\tau)={A}_{\xi}(t)$ for any $\tau \leq t$; and

(ii) $\mathcal{A}_{\xi}$ is pullback attracting, i.e.,
$$
\lim _{\tau \rightarrow-\infty} {dist}_{\mathcal H_{t}}\left(U(t, t-\tau) B, {A}_{\xi}(t)\right)=0,
$$
for any bounded subset $B \in \mathcal H_{t}(\Omega)$.
\label{def2.8-2}
\end{Definition}

\begin{Definition} {\rm(\cite{psz,zs2})}
A family ${\mathcal{A}}_{\mathcal D}=\{\mathcal{A}_{\mathcal D}(t): t \in \mathbb{R}\} \subset \Gamma\left(\mathcal{H}_{t}\right)$
is said to be a minimal time-dependent pullback $\mathcal D$-attractor for the process ${\{ U(t,\tau )\} _{t \ge \tau }}$ on $\mathcal H_{t}(\Omega)$ if the following properties hold:

(i) the set $\mathcal{A}_{\mathcal D}(t)$ is compact in $\mathcal H_{t}(\Omega)$ for any $t \in \mathbb{R}$;

(ii) ${\mathcal{A}}_{\mathcal D}$ is pullback $\mathcal D$-attracting in $\mathcal H_{t}(\Omega)$, i.e.,
$$
\lim _{\tau \rightarrow-\infty} {dist}_{\mathcal H_{t}}\left(U(t, \tau) D(\tau), \mathcal A_{\mathcal D}(t)\right)=0,
$$
for any ${\widehat D} \in {{\cal D}}$ and $t \in \mathbb{R}$;

(iii) ${\mathcal{A}}_{\mathcal D}$ is invariant, i.e., $U(t, \tau)\mathcal{A}_{\mathcal D}(\tau)={\mathcal{A}_{\mathcal D}(t)}$, for any $\tau \leq t$; and

(iv) ${\mathcal{A}}_{\mathcal D}$ is minimal, i.e., if $\widehat{C}=\{C(t): t \in \mathbb{R}\} \subset \Gamma (\mathcal{H}_{t})$ is a family of closed sets, which is pullback $\mathcal{D}$-attracting, then $\mathcal{A}_{\mathcal{D}}(t) \subset C(t)$ for any $t \in \mathbb{R}$.
\label{def2.9-2}
\end{Definition}
\begin{Remark}
The uniqueness of the minimal pullback attractor can be derived from (iv).
\label{re2.2-2}
\end{Remark}

Garc\'{i}a-Luengo, Mar\'{i}n-Rubio and Real \cite{gmr.2} proved the following lemma, which is a direct method to obtain the existence of the minimal time-dependent pullback $\mathcal D$-attractors.

\begin{Lemma}{\rm(\cite{gmr.2})}
Assume that ${\{ U(t,\tau )\} _{t \ge \tau }}$ is closed on $\mathcal H_{t}(\Omega)$, $\mathcal D$ is a universe in $\Gamma (\mathcal{H}_{t})$, $\widehat{D}_{0}=\left\{D_{0}(t): t \in \mathbb{R}\right\} \subset \Gamma (\mathcal{H}_{t})$ is a pullback $\mathcal{D}$-absorbing family for ${\{ U(t,\tau )\} _{t \ge \tau }}$ and ${\{ U(t,\tau )\} _{t \ge \tau }}$ is pullback $\widehat{D}_{0}$-asymptotically compact, then the family ${\mathcal{A}}_{\mathcal D}=\{\mathcal{A}_{\mathcal D}(t): t \in \mathbb{R}\}$ with
\[{{\cal A}_{\cal D}} = {\overline {\mathop  \bigcup \limits_{\hat D \in {\cal D}} \Lambda (t,\hat D)} ^{{{\mathcal H}_t}}}= {\overline {\mathop  \bigcup \limits_{\hat D \in \cal D} \mathop  \bigcap \limits_{s \le t} {{\overline {\mathop  \bigcup \limits_{\tau  \le s} U(t,\tau )D(\tau )} }^{{{\cal H}_t}}}} ^{{{\cal H}_t}}}\]
is the minimal pullback $\mathcal D$-attractor for the process ${\{ U(t,\tau )\} _{t \ge \tau }}$. In addition, if $\widehat{D}_{0} \in \mathcal{D}$, then
$$
\mathcal{A}_{\mathcal{D}}(t)={\mathop  \bigcap \limits_{s \le t} {{\overline {\mathop  \bigcup \limits_{\tau  \le s} U(t,\tau )D(\tau )} }^{{{\cal H}_t}}}}\subset \overline{D_{0}(t)}^{\mathcal H_{t}}
$$
for all $t \in \mathbb{R}$.
\label{th2.1-2}
\end{Lemma}

\begin{Corollary}{\rm(\cite{chm,psz})}
Let $\mathcal {D}_{F}^{\mathcal{H}_{t}}$ be the universe of fixed nonempty bounded subsets of $\mathcal{H}_{t}(\Omega)$. Namely, $\mathcal {D}_{F}^{\mathcal{H}_{t}}$ is the class of all families $\widehat{D}=\left\{D(t): t \in \mathbb{R}\right\}$, where $D(t)$ is a fixed nonempty bounded subsets of $\mathcal{H}_{t}(\Omega)$.
\label{re2.3-2}
\end{Corollary}

\begin{Corollary}{\rm(\cite{chm,psz})}\label{re2.4-2}
Under the assumptions of Lemma $\ref{th2.1-2}$, if $\mathcal {D}_{F}^{\mathcal{H}_{t}} \subset \cal D$, then the minimal time-dependent pullback  attractors ${{\cal A}_{{\cal D}_F^{{{\cal H}_t}}}}$ and $A_{\cal D}$ exist and satisfy ${{\cal A}_{{\cal D}_F^{{{\cal H}_t}}}} \subset A_{\cal D}$ for all $t \in \mathbb R$. Besides, if for some $T \in \mathbb R$ the set $\mathop  \bigcup \limits_{t \le T} {D_0}(t)$ is a bounded subset of ${\mathcal{H}_{t}(\Omega)}$, then ${{\cal A}_{{\cal D}_F^{{{\cal H}_t}}}} = A_{\cal D}$ for all ${t \le T}$.
\end{Corollary}

In order to prove the upper semicontinuity of the time-dependent pullback attractors and the global attractor, it is necessary to introduce the following lemmas.

Let $S(t): X_{t} \rightarrow X_{t}$ be a $C_{0}$-semigroup on $\left\{X_{t}\right\}_{t \in \mathbb{R}}$ and assume there exists a global attractor $A$ for $S(t)$. Then we use a non-autonomous term depending on a small parameter $\xi \in (0,\xi_{0}]$ to perturb the semigroup $S(t)$, so as to obtain a non-autonomous dynamical system driven by the process ${\{ U_{\xi}(t,\tau )\} _{t \ge \tau }}$. Moreover, for any $t \in \mathbb R$, $\tau \in \mathbb R^{+}$ and $x \in X_{t}$, assume that
\begin{equation}
\lim _{\xi \rightarrow 0} d_{X_{t}}\left(U_{\xi}(t, t-\tau) x, S(\tau) x\right)=0
\label{eq2.1-2}
\end{equation}
uniformly on any bounded set of $\left\{X_{t}\right\}_{t \in \mathbb{R}}$.

\begin{Lemma}{\rm(\cite{clr.2})}\label{th2-2}
Assume that $\rm {(\ref{eq2.1-2})}$ holds and for any small $\xi \in (0,\xi_{0}]$, there exists a pullback attractor $\mathscr{A}_{\xi}=\{{A}_{\xi}(t)\}_{t \in \mathbb R}$ and a compact set $K \subset \left\{X_{t}\right\}_{t \in \mathbb{R}}$ such that
\begin{equation}\label{eq2.2-2}
\lim _{\xi \rightarrow 0} \operatorname{dist}_{X_{t}}\left(A_{\xi}(t), K\right)=0, \quad
\end{equation}
then the upper semicontinuty of the attractors holds, that is,
\begin{equation}
\lim _{\xi \rightarrow 0} \operatorname{dist}_{X_{t}}\left(A_{\xi}(t), A\right)=0.
\label{eq2.3-2}
\end{equation}
\end{Lemma}

\begin{Lemma} {\rm(\cite{wq.2})} Assume the family ${\cal D_{\xi}}{\rm{ = }}{\{ D_{\xi}(t)\} _{t \in \mathbb R}}$ is pullback absorbing for process ${\{ U(t,\tau )\} _{t \ge \tau }}$ on $\left\{X_{t}\right\}_{t \in \mathbb{R}}$, and for any $\xi \in (0,\xi_{0}]$, the set ${\cal K_\xi } = {\{ {K_\xi }(t)\} _{t \in \mathbb R}}$ is a family of compact sets in $\left\{X_{t}\right\}_{t \in \mathbb{R}}$. Suppose ${U_\xi }( \cdot , \cdot ) = {U_{1,{\rm{ }}{\kern 1pt} \xi }}( \cdot , \cdot ) + {U_{2,{\kern 1pt} {\rm{ }}\xi }}( \cdot , \cdot ): \mathbb R \times \mathbb R \times {X_t} \to {X_t}$ such that

(i) for any $t \in \mathbb R$, $\xi \in (0,\xi_{0}]$, $x_{t-\tau} \in D_{\xi}(t-\tau)$ and $\tau>0$ the following inequality holds
\begin{equation}
\left\|U_{1,\,\xi}(t, t-\tau) x_{t-\tau}\right\|_{X_{t}} \leq \Phi(t, \tau),
\label{eq2.4-2}
\end{equation}
where $\Phi(\cdot, \cdot): \mathbb R \times \mathbb R \to \mathbb R^{+}$ and satisfies $\mathop {\lim }\limits_{\tau  \to \infty } \Phi (t,\tau ) = 0$; and

(ii) for any $t \in \mathbb R$, $\xi \in (0,\xi_{0}]$ and $T \ge 0$, the set $\mathop \bigcup \limits_{0 \le \tau  \le T} { {U}_{2,{\kern 1pt} {\rm{ }}\xi }}(t,t - \tau )D_{\xi}(t - \tau )$ is bounded and there exists a $T_{\mathcal{D_{\xi}}}(t)>0$ such that
\begin{equation}
U_{2,\, \xi}(t, t-\tau) D_{\xi}(t-\tau) \subset K_{\xi}(t)
\label{eq2.5-2}
\end{equation}
for all $\tau  \ge {T_{{{\cal D}_\xi }}}(t)$ and there exists a compact set $K \subset X_{t}$ such that
\begin{equation}
\lim _{\xi \rightarrow 0} \operatorname{dist}_{X_{t}}\left(K_{\xi}(t), K\right)=0.
\label{eq6-2}
\end{equation}
Then for each $\xi \in (0,\xi_{0}]$ there exists a pullback attractor and $\rm (\ref{eq2.2-2})$ holds.
\label{1m2.2-2}
\end{Lemma}

\section{Existence and uniqueness of solutions}
$\ \ \ \ $ Studying the attractors of an equation generally requires proving the existence and uniqueness of solutions.
Thus in this section we shall first discuss the existence and uniqueness of weak solutions to problem $(\ref{1.1-2})$.

\begin{Definition}
A weak solution to problem $(\ref{1.1-2})$ is a function $u \in C\left([\tau, T], \mathcal{H}_{t}\right(\Omega))$ for any $\tau < T$, with $u(x, \tau)=u_{\tau}$, and such that
\begin{equation}
\begin{aligned}
& \frac{d}{d t}[(u(t), \varphi)+\varepsilon(t)(\nabla u(t), \nabla \varphi)]+\left(2 a(l(u))-\varepsilon^{\prime}(t)\right)(\nabla u(t), \nabla \varphi) \\
 & = 2(f(u(t)), \varphi)+2\xi (h(t),\varphi)
\end{aligned}
\label{3.1-2}
\end{equation}
for all test functions $\varphi \in H_{0}^{1}(\Omega)$.
\label{def3.1-2}
\end{Definition}
\begin{Remark}
The equation $(\ref{3.1-2})$ should be understood in the sense of the generalized function space
$ \mathcal{D}^{\prime}(\tau,+\infty)$.
\end{Remark}
\begin{Corollary}
If $u(x,t)$ is a weak solution of problem $(\ref{1.1-2})$, then the following energy equality holds
\begin{equation}
\begin{array}{l}
\|u(t)\|^{2}+\varepsilon(t)\|\nabla u(t)\|^{2}+\int_{s}^{t}\left(2 a(l(u))-\varepsilon^{\prime}(r)\right)\|\nabla u(r)\|^{2} d r \\
=\|u(s)\|^{2}+\varepsilon(s)\|\nabla u(s)\|^{2}+2 \int_{s}^{t}(f(u(r)), u(r)) d r+2 \xi \int_{s}^{t}(h(r), u(r)) d r,
\end{array}
\label{3.2-2}
\end{equation}
for all $\tau  \le s \le t$.
\end{Corollary}
\qquad We first prove the following lemma.
\begin{Lemma} If the function $u$ is bounded in $L^{\infty}\left(\tau, T ; H_{0}^{1}(\Omega)\right)$, then $f(u)$ is bounded in $L^{q}(\tau, T ; L^{q}(\Omega))$, where $q=\frac{2(N+2)}{N \gamma}$ with $0<\gamma<\min \left\{\frac{N+2}{N-2}, 2{\rm{  +  }}\frac{4}{N}\right\}$.
\label{lem3.2-2}
\end{Lemma}
$\mathbf{Proof.}$ From the assumptions (\ref{1.7-2}), (\ref{1.8-2}), (\ref{1.10-2}) and $p \le \frac{4}{{N - 2}}$, it follows that there exists a constant satisfying $0<\gamma<\min \left\{\frac{N+2}{N-2}, 2{\rm{  +  }}\frac{4}{N}\right\}$ such that
\begin{equation}
|f(u)|\leqslant| f_{0}(u)|+| f_{1}(u)| \leqslant C\left(|u|^{\gamma}+1\right).
\label{3.3-2}
\end{equation}

Let $\gamma q=\frac{2 N+4}{N}$, then we can conclude that $q > 1$ and $2 < \gamma q < \frac{{2N}}{{N - 2}}$. Meanwhile, take $\theta {\rm{ = }}\frac{{N - 2}}{N}$, then $\theta  \in (0,1)$. Letting $\gamma q = \frac{{2N\theta }}{{N - 2}} + 2(1 - \theta )$ and using (\ref{3.3-2}) and the H$\ddot{o}$lder inequality, we obtain
\begin{equation}
\begin{aligned}
\|f(u)\|_{L^{q}(\Omega)}^{q} & \leqslant C+C \int_{\Omega}|u|^{\mid \gamma q} d x \\
& \leqslant C+C \int_{\Omega}|u|^{\frac{2 N}{N-2}} d x \\
&=C+C \int_{\Omega}|u|^{\frac{2 N \theta}{N-2}+2(1-\theta)} d x \\
& \leqslant C+C(\int_{\Omega}|u|^{\frac{2 N}{N-2}} d x)^{\theta}(\int_{\Omega}|u|^{2} d x)^{1-\theta}.
\end{aligned}
\label{3.4-2}
\end{equation}

Then from (\ref{3.4-2}) and since when $N \ge 3$ the embedding $H_{0}^{1}(\Omega) \subset L^{\frac{2 N}{N-2}}(\Omega)$ is continuous (see \cite{Evans}), the following inequalities hold
\begin{equation}
\begin{gathered}
\|f(u)\|_{L^{q}(\Omega)}^{q} \leqslant C+C\|u\|_{L^{\frac{2 N}{N-2}}(\Omega)}^{\frac{2 N \theta}{N-2}}\|u\|_{L^{2}(\Omega)}^{2(1-\theta)} \\
\quad\leqslant C+C\|u\|_{L^{\frac{2 N}{N-2}}(\Omega)}^{\frac{2 N \theta}{N-2}} \\
\quad=C+C\|u\|_{L^{\frac{2 N}{N-2}}(\Omega)}^{2} \\
\qquad \leqslant C+C\|\nabla u\|_{L^{\frac{2 N}{N-2}}(\Omega)}^{2}.
\end{gathered}
\label{3.5-2}
\end{equation}

The proof is complete by (\ref{3.5-2}). $\hfill$$\Box$

Now we prove the existence and uniqueness of the solutions to problem $(\ref{1.1-2})$.
\begin{Theorem}
Assume that $a(\, \cdot \,)$ is a local Lipschitz continuous function and satisfies $(\ref{1.4-2})$, $l(\, \cdot \,)$ is given in $(\ref{Lg-2})$, $f \in C(\mathbb{R}, \mathbb{R})$ and satisfies $(\ref{1.5-2})-(\ref{1.10-2})$, $h \in L_{l o c}^{2}(\mathbb{R} ; L^{2}(\Omega))$ and the initial value ${u_\tau } \in {{\cal H}_t}(\Omega )$, then for any   $\tau \in \mathbb{R}$ and $t \ge \tau $, there exists a weak solution to problem $(\ref{1.1-2})$. Moreover, the solution $u$ depends continuously on its initial value.
\label{th3.1-2}
\end{Theorem}
$\mathbf{Proof.}$ Consider the approximate solution $u_{k}\left(t, \tau ; u_{\tau}\right)$ $=\sum\limits_{j=1}^{k} r_{k, j}(t) \omega_{j}(x)$, where $j, k \in {\mathbb N^+}$, $\left\{ {{\omega_j}} \right\}_{j = 1}^\infty $ is a basis of $H^{2}(\Omega) \cap H_{0}^{1}(\Omega)$ and orthonormal in $L^{2}(\Omega)$. The Faedo-Galerkin method needs to find an approximate sequence $\{ {u_k}\} $ so that the following approximate system holds:
\begin{equation}
\left\{ {\begin{array}{*{20}{l}}
{\frac{d}{{dt}}[({u_k}(t),{\omega _j}) + \varepsilon (t)(\nabla {u_k}(t),\nabla {\omega _j})] + ( {2a(l({u_k})) - {\varepsilon ^\prime }(t)} )(\nabla {u_k}(t),\nabla {\omega _j})}\\
{ = 2(f({u_k}(t)),{\omega _j}) + 2\xi\left({h(t),{\omega _j}} \right), \quad \forall\, t \in [\tau , + \infty ),}\\
{(u_{k}(\tau),\omega_{j})=(u_{\tau},\omega_{j}), \quad j=1, 2, \cdots, k.}
\end{array}} \right.
\label{3.6-2}
\end{equation}

\textbf{\textbf{\emph{Step}}\,\emph{1:}\,\emph{(A priori estimate for $\boldsymbol{u_k}$)}}
Multiplying $(\ref{3.6-2})_{1}$ by the test function ${\gamma _{k,{\rm{ }}j}}(t)$ and then summing $j$ from $1$ to $k$, we obtain
\begin{equation}
\begin{aligned}
&\frac{d}{d t}(\left\|u_{k}(t)\right\|^{2}+\varepsilon(t)\left\|\nabla u_{k}(t)\right\|^{2})+(2a(l(u))-\varepsilon^{\prime}(t))\left\|\nabla u_{k}(t)\right\|^{2}\\
&= 2\left(f\left(u_{k}(t)\right), u_{k}(t)\right)+2\xi\left(h(t), u_{k}(t)\right).
\end{aligned}
\label{3.7-2}
\end{equation}

From (\ref{1.6-2}) there exists some constants $0<\tilde \eta < \frac{3}{4}m{\lambda _1} $ and $C_{0}>0$ such that
\begin{equation}
2\left(f\left(u_{k}(t)\right), u_{k}(t)\right) \leqslant (\frac{3}{2} m \lambda_{1}-2 \tilde{\eta} )\left\|u_{k}(t)\right\|^{2}+ 2 C_{0}.
\label{3.8-2}
\end{equation}

Using the Young and the Poincar$\acute{e}$ inequalities, it follows that
\begin{equation}
2 \xi\left(h(x, t), u_{k}(t)\right) \leq \frac{2 \xi^{2}}{m \lambda_{1}}\|h(x, t)\|^{2}+\frac{m}{2}\left\|\nabla u_{k}(t)\right\|^{2}.
\label{3.9-2}
\end{equation}

Substituting (\ref{3.8-2}) and (\ref{3.9-2}) into (\ref{3.7-2}), and then by (\ref{1.4-2}) and the Poincar$\acute{e}$ inequality, we can derive
\begin{equation}
\frac{d}{d t}(\|u_{k}(t)\|^{2}+\varepsilon(t)\|u_{k}(t)\|^{2})+\left(2 \tilde{\eta}-\varepsilon^{\prime}(t)\right)\left\|\nabla u_{k}(t)\right\|^{2} \leqslant \frac{2 \xi^{2}}{m \lambda_{1}}\|h(x, t)\|^{2}+2 C_{0}.
\label{3.10-2}
\end{equation}

Integrating (\ref{3.10-2}) from $\tau$ to $t$, we deduce
\begin{equation}
\begin{aligned}
&\left\|u_{k}(t)\right\|^{2}+\varepsilon(t)\left\|\nabla u_{k}(t)\right\|^{2}+\int_{\tau}^{t}(2 \tilde{\eta}-\varepsilon^{\prime}(s))\left\|\nabla u_{k}(s)\right\|^{2} d s \\
&\leqslant\left\|u_{k}(\tau)\right\|^{2}+\varepsilon(\tau)\left\|\nabla u_{k}(\tau)\right\|^{2}+\frac{2 \xi^{2}}{m \lambda_{1}} \int_{\tau}^{t}\|h(x, s)\|^{2} d s+2 C_{0}(t-\tau) .
\end{aligned}
\label{3.11-2}
\end{equation}

From (\ref{3.11-2}) and notice that $\varepsilon(t)$ is a decreasing bounded function, we obtain that for any $T > t$, $\{ {u_k}\} $ is bounded in $L^{\infty}(\tau, T ; L^{2}(\Omega)) \cap L^{\infty}(\tau, T ; H_{0}^{1}(\Omega)) \cap L^{2}(\tau, T ; H_{0}^{1}(\Omega))$, then we deduce
\begin{equation}
\left\{u_{k}\right\} \text { is bounded in } L^{\infty}(\tau, T ; \mathcal{H}_{t}(\Omega)).
\label{3.12-2}
\end{equation}

Multiplying the approximate system $(\ref{3.6-2})_{1}$ by $\gamma_{k, j}^{\prime}(t)$ and summing $j$ from 1 to $k$, then by (\ref{1.4-2}), we obtain
\begin{equation}
\left\|\left(u_{k}(t)\right)_{t}\right\|^{2}+\varepsilon(t)\left\|\nabla\left(u_{k}(t)\right) _{t}\right\|^{2}+\frac{m}{2} \frac{d}{d t}\left\|\nabla u_{k}(t)\right\|^{2} \leqslant\left(f\left(u_{k}(t)\right), u_{k}(t)\right)+\xi\left(h(x, t), u_{k}(t)\right).
\label{3.13-2}
\end{equation}

Using the Young inequality, we obtain
\begin{equation}
\frac{1}{2}\left\|\left(u_{k}(t)\right)_{t}\right\|^{2}+\varepsilon(t)\left\|\nabla\left(u_{k}(t)\right)_{t}\right\|^{2}+\frac{m}{2} \frac{d}{d t}\left\|\nabla u_{k}(t)\right\|^{2} \leqslant \frac{1}{4}\left\|f\left(u_{k}(t)\right)\right\|^{2}+\frac{\xi^{2}}{4}\left\|h\left(x, t\right)\right\|^{2}.
\label{3.14-2}
\end{equation}

Integrating (\ref{3.14-2}) from $\tau$ to $t$, we deduce
\begin{equation}
\begin{aligned}
&\int_{\tau}^{t}(\frac{1}{2}\|(u_{k}(s))_{s}\|^{2}+\varepsilon(s)\|\nabla(u_{k}(s))_{s}\|^{2}) d s+\frac{m}{2}\|\nabla u_{k}(t)\|^{2} \\
&\leqslant \frac{m}{2}\left\|\nabla u_{k}(\tau)\right\|^{2}+\frac{1}{4} \int_{\tau}^{t}\left\|f\left(u_{k}(s)\right)\right\|^{2} d s+\frac{\xi^{2}}{4} \int_{\tau}^{t}\|h(x, s)\|^{2} d s.
\end{aligned}
\label{3.15-2}
\end{equation}

By Lemma \ref{lem3.2-2}, it follows that
\begin{equation}
\left\{f(u_{k})\right\} \text { is bounded in } L^{q}(\tau, T ; H_0^1(\Omega)).
\label{3.16-2}
\end{equation}

Then from (\ref{3.15-2}) and (\ref{3.16-2}), through similar calculations and estimations to (\ref{3.12-2}), we deduce
\begin{equation}
\left\{\partial _{t}u_{k}\right\} \text { is bounded in } L^{\infty}(\tau, T ; \mathcal H_{t}(\Omega)).
\label{3.17-2}
\end{equation}

From (\ref{3.12-2}), (\ref{3.16-2}), (\ref{3.17-2}), the compactness arguments and the Aubin-Lions lemma (see \cite{Lions}), we derive that there exists a subset of $\left\{u_{k}\right\}$ (still marked as $\left\{u_{k}\right\}$), $u \in L^{\infty}\left(\tau, T ; \mathcal{H}_{t}(\Omega)\right) \cap L^{2}(\tau, T ; H_{0}^{1}(\Omega))$ and $\partial_{t} u \in L^{\infty}\left(\tau, T ; \mathcal{H}_{t}(\Omega)\right)$ such that
\begin{equation}
\quad \quad u_{k} \rightharpoonup u \quad \text { weakly-star in } L^{\infty}(\tau, T ;\mathcal H_{t}(\Omega));
\label{3.18-2}
\end{equation}
\begin{equation}
u_{k} \rightharpoonup u \quad \text { weakly in } L^{2}(\tau, T ; H_0^1(\Omega ));
\label{3.19-2}
\end{equation}
\begin{equation}
\quad f\left(u_{k}\right) \rightharpoonup f\left(u\right) \quad \text { weakly in } L^{q}(\tau, T ; L^{q}(\Omega));
\label{3.20-2}
\end{equation}
\begin{equation}
\quad\quad\quad\quad a\left(l(u_{k})\right) u_{k} \rightharpoonup a\left(l(u)\right) u \quad \text { weakly in } L^{2}(\tau, T ; H_{0}^{1}(\Omega));
\label{3.21-2}
\end{equation}
\begin{equation}
\quad {\partial _t}{u_k} \rightharpoonup {\partial _t}u \quad \text { weakly in } L^{2}(\tau, T ;\mathcal H_{t}(\Omega));
\label{3.22-2}
\end{equation}
\begin{equation}
u_{k} \rightarrow u \quad a. e.\quad(x, t) \in \Omega \times[\tau,+\infty).
\label{3.23-2}
\end{equation}

\textbf{\textbf{\emph{Step}}\,\emph{2:}\,\emph{(Verify the continuity of $\boldsymbol u$)}} Let $\tilde u=u_{k}-u$, then $\tilde u$ satisfies
\begin{equation}
\begin{aligned}
&\frac{d}{d t}\|\tilde{u}\|_{\mathcal H_{t}}^{2}=\frac{d}{d t}\left[\|\tilde{u}\|^{2}+\varepsilon(t)\|\nabla \tilde{u}\|^{2}\right]+2(a(l(u_{k}))-\varepsilon^{\prime}(t))\|\nabla \tilde{u}\|^{2} \\
&=2\left(a(l(u))-a\left(l\left(u_{k}\right)\right)\right)(\nabla u, \nabla \tilde{u})+2\left(f\left(u_{k}\right)-f(u), \tilde{u}\right).
\end{aligned}
\label{3.24-2}
\end{equation}

Noting that the function $a(\, \cdot \,)$ is local Lipschitz continuous and by the Young and the Cauchy inequalities, we can obtain
\begin{equation}
2\left(a(l(u))-a\left(l\left(u_{k}\right)\right)\right)(\nabla u, \nabla \tilde{u}) \leqslant 2 m\left\|\nabla u_{k}-\nabla u\right\|^{2}+\frac{(L a(R))^{2}\|l\|^{2}\|\nabla u\|^{2}\left\|u_{k}-u\right\|^{2}}{2 m},
\label{3.25-2}
\end{equation}
where $L_{a}(R)$ is the Lipschitz constant of the function $a(\, \cdot \,)$ in $[-R, R]$.

By (\ref{1.5-2}), the H$\ddot{o}$lder, the Poincar$\acute{e}$ inequalities, definition of the space $\mathcal H_{t}(\Omega)$ and notice that $p \le \frac{4}{{N - 2}}$, then for any $0<\alpha<\min \left\{1, \frac{4-(N-2) p}{2}\right\}$, we obtain
\begin{equation}
\begin{aligned}
2\left(f\left(u_{k}\right)-f(u), \tilde{u}\right) & \leq C \int_{\Omega}\left(1+\left|u_{k}\right|^{p}+|u|^{p}\right)|\tilde{u}|^{2} d x \\
& \leq C(\int_{\Omega}(1+|u_{k}|^{p}+|u|^{p})^{\frac{2 N}{(N-2) p}} d x)^{\frac{(N-2) p}{2 N}}(\int_{\Omega}|\tilde{u}|^{\frac{2 N}{N-2}} d x)^{\frac{N-2}{2 N}} \\
& \times(\int_{\Omega}|\tilde{u}|^{\frac{2 N}{2N-(N-2)(p+1)}} d x)^{\frac{2N-(N-2)(p+1)}{2 N}} \\
& \leq C(1+\|A^{\frac{1}{2}} u_{k}\|^{p}+\|A^{\frac{1}{2}} u\|^{p})\|A^{\frac{1}{2}} \tilde{u}\|\|A^{\frac{1-\alpha}{2}} \tilde{u}\| \\
&\leq C(1+\|A^{\frac{1}{2}} u_{k}\|^{p}+\|A^{\frac{1}{2}} u\|^{p})(\|\tilde{u}\|^{2}+\varepsilon(t)\|\nabla \tilde{u}\|^{2}) \\
& \leq C (1+ \|A^{\frac{1}{2}} u_{k} \|^{\frac{4}{N-2}}+ \|A^{\frac{1}{2}} u \|^{\frac{4}{N-2}} )\|\tilde{u}\|_{\mathcal{H}_{t}}^{2}.
\end{aligned}
\label{3.26-2}
\end{equation}

Substituting (\ref{3.25-2}) and (\ref{3.26-2}) into (\ref{3.24-2}), from (\ref{1.4-2}), $\varepsilon(t)$ is a decreasing function and using the Sobolev embedding theorem (see \cite{Evans}), we can derive
\begin{equation}
\frac{d}{d t}\|\tilde{u}\|_{\mathcal H_{t}}^{2}=\frac{d}{d t}\left(\|\tilde{u}\|^{2}+\varepsilon(t)\|\nabla \tilde{u}\|^{2}\right) \leqslant C\left(\|\tilde{u}\|^{2}+\varepsilon(t)\|\nabla \tilde{u}\|^{2}\right).
\label{3.27-2}
\end{equation}

Applying the generalized Gronwall lemma (see \cite{Q1,Q2,Q3}) to $(\ref{3.27-2})$, we can conclude
\begin{equation}
\|\tilde{u}\|_{\mathcal H_{t}}^{2}=\|\tilde{u}(t)\|^{2}+\varepsilon(t)\|\nabla \tilde{u}(t)\|^{2} \leqslant e^{C(t-\tau)}(\left\|\tilde{u}_{\tau}\|^{2}+\varepsilon(\tau)\|\tilde{u}_{\tau}\|^{2}\right).
\label{3.28-2}
\end{equation}

\textbf{\textbf{\emph{Step}}\,\emph{3:}\,\emph{(Verify the initial value $\boldsymbol {u_{\tau}}$)}} Choosing a suitable test function $\varphi \in C^{1}([\tau, T] ; H_{0}^{1}(\Omega))$ with $\varphi(T)=0$, then it follows that
\begin{equation}
\begin{aligned}
&\int_{\tau}^{T}-\left(u, \varphi^{\prime}\right) d s+\int_{\tau}^{T} \int_{\Omega} \varepsilon(s) \nabla u_{s} \nabla \varphi d x d s-\int_{\tau}^{T} \int_{\Omega} a(l(u)) (\Delta u) \varphi d x d s\\
&-\int_{\tau}^{T} \int_{\Omega}(f(u)+\xi h(x, s))\varphi d x d s =(u(\tau), \varphi(\tau)) .
\end{aligned}
\label{3.29-2}
\end{equation}

In the same way as in the Faedo-Galerkin approximations, we conclude
\begin{equation}
\begin{aligned}
&\int_{\tau}^{T}- (u_{k}, \varphi^{\prime} ) d s+\int_{\tau}^{T} \int_{\Omega} \varepsilon(s) \nabla (u_{k} )_{s} \nabla \varphi d x d s-\int_{\tau}^{T} \int_{\Omega} a(l(u_{k})) (\Delta u_{k}) \varphi d x d s \\
&-\int_{\tau}^{T} (f (u_{k})+\xi h(x, s))\varphi d s =\left(u_{k}(\tau), \varphi(\tau)\right).
\end{aligned}
\label{3.30-2}
\end{equation}

Taking limits as $k \to \infty $ in (\ref{3.30-2}) and since $u_{k}(\tau) \rightarrow u_{\tau}$, we obtain
\begin{equation}
\begin{aligned}
&\int_{\tau}^{T}- (u, \varphi^{\prime} ) d s+\int_{\tau}^{T} \int_{\Omega} \varepsilon(s) \nabla u_{s} \nabla \varphi d x d s-\int_{\tau}^{T} \int_{\Omega} a(l(u)) (\Delta u) \varphi d x d s \\
&-\int_{\tau}^{T} (f (u)+\xi h(x, s))\varphi d s =\left(u_{\tau}, \varphi(\tau)\right).
\end{aligned}
\label{3.31-2}
\end{equation}

Then $u(\tau)=u_{\tau}$ directly holds.

From the above estimations, we can obtain that $u$ is a weak solution of problem $(\ref{1.1-2})$. $\hfill$$\Box$

\begin{Theorem}
Under the assumptions of Theorem ${\ref{th3.1-2}}$, if the weak solution of problem $(\ref{1.1-2})$ exists, then it is a unique solution.
\label{th3.2-2}
\end{Theorem}
$\mathbf{Proof.}$ Assuming that $u^{1}$ and $u^{2}$ are two solutions corresponding to the initial values $u_{\tau}^{1}$ and $u_{\tau}^{2}$, respectively, and satisfying
\begin{equation}
\left\{\begin{array}{ll}
u^{1}_{t}-\varepsilon(t) \Delta u^{1}_{t}-a(l(u^{1})) \Delta u^{1}=f(u^{1})+\xi h(x,t) & \text { in } \Omega \times(\tau, \infty), \\
u^{1}=0 & \text { on } \partial \Omega\times(\tau, \infty), \\
u^{1}(x, \tau)=u_{\tau}^{1}(x),  &\,\, x \in \Omega,
\end{array}\right.
\label{3.32-2}
\end{equation}
and
\begin{equation}
\left\{\begin{array}{ll}
u^{2}_{t}-\varepsilon(t) \Delta u^{2}_{t}-a(l(u^{2})) \Delta u^{2}=f(u^{2})+\xi h(x,t) & \text { in } \Omega \times(\tau, \infty), \\
u^{2}=0 & \text { on } \partial \Omega\times(\tau, \infty), \\
u^{2}(x, \tau)=u_{\tau}^{2}(x),  &\,\, x \in \Omega.
\end{array}\right.
\label{3.33-2}
\end{equation}

Subtracting  (\ref{3.33-2}) from (\ref{3.32-2}) and taking $u=u^{1}-u^{2}$ as the test function of the resulting equation, we derive
\begin{equation}
\begin{aligned}
&\frac{d}{d t}[\|u\|^{2}+\varepsilon(t)\|\nabla u\|^{2}]+2(a(l(u^{1}))-\varepsilon^{\prime}(t))\|\nabla u\|^{2} \\
&=2(a(l(u^{2}))-a(l(u^{1})))(\nabla u^{2}, \nabla u)+2(f(u^{1})-f(u^{2}), u).
\end{aligned}
\label{3.34-2}
\end{equation}

Performing similar calculations to the proof of Theorem ${\ref{th3.1-2}}$ can obtain the following inequality
\begin{equation}
\|{u}\|_{\mathcal H_{t}}^{2}=\|{u}(t)\|^{2}+\varepsilon(t)\|\nabla {u}(t)\|^{2} \leqslant e^{C(t-\tau)}(\|{u}_{\tau}\|^{2}+\varepsilon(\tau)\|{u}_{\tau}\|^{2}).
\label{3.35-2}
\end{equation}

Consequently, the uniqueness of the solution follows readily.$\hfill$$\Box$

\begin{Corollary}
Thanks to Theorems $\rm{\ref{th3.1-2}}$ and $\rm{\ref{th3.2-2}}$, problem $(\ref{1.1-2})$ has a continuous process
$$
U(t, \tau):\mathcal{ H}_{t}(\Omega) \rightarrow \mathcal{ H}_{t}(\Omega)
$$
with $U(t, \tau) u_{\tau}$ being the unique weak solution of $\rm({\ref{1.1-2}})$ respects to initial datum $u_{\tau}$.
\label{re3.1-2}
\end{Corollary}

\section{Existence of the minimal time-dependent pullback $\mathcal D$-attractors}
\ \ \ \ In this section, we will verify the existence of the minimal time-dependent pullback $\mathcal D$-attractors for the process  $\left\{ {U(t,\tau ){\} _{t \ge \tau }}} \right.$ in $\mathcal H_{t}(\Omega)$. To prove it, we will check the four properties mentioned in Definition \ref{def2.9-2}, therefore we first estimate the following lemma.
\begin{Lemma}\label{lem4.1-2}
Under the assumptions of Theorems ${\ref{th3.1-2}}$ and ${\ref{th3.2-2}}$, then for any $t \ge \tau $, the solution of problem $(\ref{1.1-2})$ satisfies
\begin{equation}
\|{u(t)}\|_{\mathcal H_{t}}^{2}=\|u(t)\|^{2}+\varepsilon(t)\|\nabla u(t)\|^{2} \leq e^{-\sigma \tau}\left\|u_{t-\tau}\right\|_{\mathcal H_{t}}^{2}+\frac{\xi e^{-\sigma t}}{\eta} \int_{{\rm{ - }}\infty }^{t} e^{\sigma s}\|h(x, s)\|^{2} d s +\frac{2C_{1}}{\sigma},
\label{4.1-2}
\end{equation}
where
$0<\sigma<\delta_{1}<\min \left\{\eta,\frac{{{\rm{ - }}{\varepsilon^{\prime}}(t)}}{{\varepsilon (t)}} \right\}$
with $0 < \eta  < m{\lambda _1}$.
\end{Lemma}
$\mathbf{Proof.}$ From $(\ref{3.2-2})$, it easily follows that the weak solution $u$ satisfies
\begin{equation}
\begin{aligned}
&\frac{d}{d t}\left[\|u(t)\|^{2}+\varepsilon(t)\|\nabla u(t)\|^{2}\right]+\left(2 a(l(u(t)))-\varepsilon^{\prime}(t)\right)\|\nabla u(t)\|^{2} \\
&=2(f(u(t)), u(t))+2\xi(h(x, t), u(t)).
\label{4.2-2}
\end{aligned}
\end{equation}

By (\ref{1.4-2}), we obtain
\begin{equation}
2 m\|\nabla  u(t)\|^{2} \leqslant 2 a(l(u))\|\nabla u(t)\|^{2}.
\label{4.3-2}
\end{equation}

Besides, we conclude from (\ref{1.6-2}) there exists $0<\eta<m \lambda_{1}$ such that
\begin{equation}
(f(u(t)), u(t)) \leqslant\left(m \lambda_{1}-\eta\right)\|u(t)\|^{2}+C_{1}.
\label{4.4-2}
\end{equation}

Using the Young and the Cauchy inequalities, we can derive
\begin{equation}
2 \xi (h(x, t), u(t)) \leqslant \frac{\xi}{\eta}\|h(x, t)\|^{2}+\xi \eta\|u(t)\|^{2}.
\label{4.5-2}
\end{equation}

Inserting $(\ref{4.3-2})-(\ref{4.5-2})$ into (\ref{4.2-2}), we obtain
\begin{equation}
\begin{aligned}
& \frac{d}{d t}\left[\|u(t)\|^{2}+\varepsilon(t)\|\nabla u(t)\|^{2}\right]+(2 m-\varepsilon^{\prime}(t))\|\nabla u(t)\|^{2} \\
& \leqslant 2\left(m \lambda_{1}-\eta\right)\|u(t)\|^{2}+2 C_{1}+\frac{\xi}{\eta}\|h(x,t)\|^{2}+\xi \eta\|u(t)\|^{2}.
\end{aligned}
\label{4.6-2}
\end{equation}

Then by the Poincar$\acute{e}$ inequality, we can derive
\begin{equation}
\frac{d}{d t}\left[\|u(t)\|^{2}+\varepsilon(t)\|\nabla u(t)\|^{2}\right]+(2 \eta-\xi \eta)\|u(t)\|^{2}-\varepsilon^{\prime}(t)\|\nabla u(t)\|^{2} \leqslant 2 C_{1}+\frac{\xi}{\eta}\|h(x,t)\|^{2}.
\label{4.7-2}
\end{equation}

Taking $\xi$ is so small that $(2 \eta-\xi \eta) \| u(t)\|^{2}>\eta\| u(t) \|^{2}$, which substituted into (\ref{4.7-2}) leads to
$$
\frac{d}{d t}[\|u(t)\|^{2}+\varepsilon(t) \|\left. \nabla u(t)\right|^{2}]+\eta\|u(t)\|^{2}-\varepsilon^{\prime}(t)\|\nabla u(t)\|^{2} \leqslant 2 C_{1}+\frac{\xi}{\eta} \| h(x,t)\|^{2}.
$$

Taking $0<\sigma<\delta_{1}<\min \left\{\eta, \frac{-\varepsilon^{\prime}(t)}{\varepsilon(t)}\right\}$ with $0 < \eta  < m{\lambda _1}$, then we can derive
\begin{equation}
\frac{d}{d t}[\|u(t)\|^{2}+\varepsilon(t) \|\left. \nabla u(t)\right||^{2}]+\delta_{1}(\|u(t)\|^{2}+\varepsilon(t){\|\nabla u(t)\|}^{2}) \leqslant 2 C_{1}+\frac{\xi}{\eta} \| h(x,t)\|^{2}.
\label{4.8-2}
\end{equation}

Then by a simple calculation, we can conclude
\begin{equation}
\begin{aligned}
& \frac{d}{d t}(e^{\sigma t}(\|u(t)\|^{2}+\varepsilon (t)\|\nabla u(t)\|^{2}))+(\delta_{1}-\sigma) e^{\sigma t}(\| u(t)\|^{2}+\varepsilon(t)\| \nabla u(t) \|^{2}) \\
& \leqslant 2 C_{1}e^{\sigma t}+\frac{\xi}{\eta} e^{\sigma t}\|h(x,t)\|^{2} .
\end{aligned}
\label{4.9-2}
\end{equation}

Integrating (\ref{4.9-2}) from $t-\tau$ to $t$, we deduce
\begin{equation}
\begin{aligned}
&\|u(t)\|^{2}+\varepsilon(t)\|\nabla u( t) \|^{2}+\left(\delta_{1}-\sigma\right) e^{-\sigma t} \int_{t-\tau}^{t} e^{\sigma s}\left(\|u(s)\|^{2}+\varepsilon(s)\|\nabla u(s)\|^{2}\right) d s \\
&\leqslant e^{-\sigma \tau}(\|u_{t-\tau}\|^{2}+\varepsilon(t-\tau) \| \nabla u_{t-\tau }\|^{2})+\frac{\xi e^{-\sigma t}}{\eta} \int_{-\infty}^{t} e^{\sigma s}\|h(x, s)\|^{2} d s+\frac{2 C_{1}}{\sigma}.
\end{aligned}
\label{4.10-2}
\end{equation}

From the definition of $\mathcal H_{t}(\Omega)$ and (\ref{1.3-2}), (\ref{4.1-2}) follows directly. $\hfill$$\Box$

\begin{Definition} (Tempered universe)\label{def4.1-2}
 For each $\sigma>0$, let $\mathcal{D}_{\sigma}^{\mathcal H_{t}}$ be the class of all families of nonempty subsets $\widehat{D}=\{D(t): t \in \mathbb{R}\} \subset \Gamma\left(\mathcal H_{t}\right)$ such that
$$
\lim _{\tau \rightarrow-\infty}\left(e^{\sigma \tau} \sup _{v \in D(\tau)}{\|v\|}^{2}_{\mathcal H_{t}}\right)=0.
$$
\end{Definition}
\begin{Remark}
The universe $\mathcal{D}_{F}^{L^{2}} \subset \mathcal{D}_{\sigma}^{\mathcal H_{t}}$ and $\mathcal{D}_{\sigma}^{\mathcal H_{t}}$ is inclusion-closed, which means that if $\widehat{D} \in \mathcal{D}_{\sigma}^{\mathcal H_{t}}$ and $\widehat{D}^{\prime}=\left\{D^{\prime}(t): t \in \mathbb{R}\right\} \subset \Gamma(\mathcal H_{t})$ satisfies that $D^{\prime}(t) \subset D(t)$ for all $t \in \mathbb{R}$, then $\widehat{D}^{\prime} \in\mathcal{D}_{\sigma}^{\mathcal H_{t}}$.
\label{re4.1-2}
\end{Remark}

Based on the above results, adding some suitable growth conditions to the function $h$ of problem (\ref{1.1-2}), then we obtain the existence of the $\mathcal{D}_{\sigma}^{\mathcal H_{t}}$-absorbing family of process ${\{ U(t,\tau )\} _{t \ge \tau }}$ on $\mathcal H_{t}(\Omega)$.
\begin{Lemma}
Under the assumptions of Theorems $\ref{th3.1-2}$ and ${\ref{th3.2-2}}$, if $h(x, t)$ also satisfies
\begin{equation}
\int_{-\infty}^{0} e^{\sigma s}\|h(x,s)\|^{2} d s<+\infty,
\label{4.11-2}
\end{equation}
for some $0<\sigma<\delta_{1}<\min \left\{\eta,\frac{{{\rm{ - }}{\varepsilon^{\prime}}(t)}}{{\varepsilon (t)}} \right\}$
with $0 < \eta  < m{\lambda _1}$. Then the family $\widehat{D}_{0}=\left\{D_{0}(t): t \in \mathbb{R}\right\}$ with $D_{0}(t)=\bar{B}_{\mathcal H_{t}}\left(0, {\rho}_{\xi}(t)\right)$, the closed ball in $\mathcal{H}_{t}(\Omega)$ of centre zero and radius $\rho_{\xi}(t)$, where
\begin{equation}
\rho_{\xi}(t)=C_{2}\left(\xi e^{-\sigma t} \int_{-\infty}^{t} e^{\sigma s}\|h(x, s)\|^{2} d s+1\right)
\label{4.12-2}
\end{equation}
and ${C_2} = \max \{ \frac{2}{\eta },\frac{{4{C_1}}}{\sigma }\} $ is pullback $\mathcal{D}_{\sigma}^{\mathcal H_{t}}$-absorbing family of process ${\{ U(t,\tau )\} _{t \ge \tau }}$ on $\mathcal H_{t}(\Omega)$. Moreover, $\widehat{D}_{0} \in \mathcal{D}_{\sigma}^{\mathcal H_{t}}$.
\label{lem4.2-2}
\end{Lemma}
$\mathbf{Proof.}$  From the prove of Lemma \ref{lem4.1-2} and Definition \ref{def4.1-2}, we can derive that Lemma \ref{lem4.2-2} follows directly. $\hfill$$\Box$

To prove the existence of the minimal time-dependent pullback $\mathcal{D}_{\sigma}^{\mathcal H_{t}}$-attractors for the process $\left\{ {U(t,\tau ){\} _{t \ge \tau }}} \right.$, we shall check the compactness of  process ${\{ U(t,\tau )\} _{t \ge \tau }}$ on $\mathcal H_{t}(\Omega)$.
\begin{Lemma}
Under the assumptions of Lemma $\ref{lem4.2-2}$, for any $t \in \mathbb R$ and $\widehat{D} \in \mathcal{D}_{\sigma}^{\mathcal H_{t}}$, there exists $\tau_{1}(t,\widehat{D})<t-2 \text { such that for any } \tau \leq \tau_{1}(t,\widehat{D}) \text { and any } u_{\tau} \in D(\tau)$, the following inequalities hold:
\begin{equation}
{\left\|u\left(r ; \tau, u_{\tau}\right)|\right|}^{2}_{\mathcal H_{t}} \leq \rho_{1}(t), \quad \forall \, r \in[t-2, t],
\label{4.13-2}
\end{equation}
\begin{equation}
\int_{r-1}^{r}\left\|\nabla u\left(s ; \tau, u_{\tau}\right)\right\|^{2} d s \leq \rho_{2}(t), \quad \forall \, r \in[t-1, t],
\label{4.14-2}
\end{equation}
where
$$
\rho_{1}(t)=C_{2}\left(1+\xi e^{-\sigma (t-\tau)} \int_{-\infty}^{t} e^{\sigma s}\|h(x, s)\|^{2} d s\right),
$$
$$
\rho_{2}(t)=\frac{1}{2\tilde{\eta}+L}\left(\rho_{1}(t)+\frac{{2{\xi ^2}}}{{m{\lambda _1}}}\max _{r \in[t-1, t]} \int_{r-1}^{r}\|h(x,s)\|^{2} d s+2C_{0}\right),
$$
$\tilde{\eta}$, $C_{0}$ and $C_{2}$ are the same as they are in Theorem $\ref{th3.1-2}$ and Lemma $\ref{lem4.2-2}$, respectively.
\label{lem4.3-2}
\end{Lemma}
$\mathbf{Proof.}$ Let $\tau \leq \tau_{1}(t,\widehat{D})<t-2$, then from Definition $\ref{def4.1-2}$ and similar calculations to Lemma $\ref{lem4.1-2}$, we can derive that $(\ref{4.13-2})$ holds.
Then through similar calculations to $(\ref{3.11-2})$, we can obtain
\begin{equation}
\frac{d}{d s}(\|u_{k}(s)\|^{2}+\varepsilon(s)\left\|\nabla u_{k}(s)\|^{2}\right)+\left(2 \tilde{\eta}-\varepsilon^{\prime}(t)\right)\left\|\nabla u_{k}(s)\right\|^{2} \leqslant \frac{2 \xi^{2}}{m \lambda_{1}}\|h(x, s)\|^{2}+2 C_{0}.
\label{4.15-2}
\end{equation}

Integrating (\ref{4.15-2}) from $r-1$ to $r$, we deduce
\begin{equation}
\begin{aligned}
&\|u_{k}(r)\|^{2}+\varepsilon (r)\|\nabla u_{k}(r)\|^{2}+\int_{r-1}^{r}(2 \tilde{\eta}-\varepsilon^{\prime}(s))\|\nabla u_{k}(s)\|^{2} d s. \\
&\leqslant \|u_{k}(r-1) \|^{2}+\varepsilon(r-1)\|\nabla u_{k}(r-1)\|^{2}+\frac{2 \xi^{2}}{m\lambda_{1}} \int_{r-1}^{r}\|h(x, s)\|^{2} d s+2 C_{0}.
\end{aligned}
\label{4.16-2}
\end{equation}

By (\ref{1.3-2}), we obtain
\begin{equation}
\int_{r-1}^{r}(2 \tilde{\eta}-\varepsilon^{\prime}(s))\left\|\nabla u_{k}(s)\right\|^{2} d s \leqslant(2 \tilde{\eta}+L) \int_{r-1}^{r}\left\|\nabla u_{k}(s)\right\|^{2} d s.
\label{4.17-2}
\end{equation}

Then from (\ref{4.13-2}), (\ref{4.16-2}) and (\ref{4.17-2}), we derive that for any $k \ge 1$
\begin{equation}
\int_{r-1}^{r}\left\|\nabla u_{k}(s)\right\|^{2} d s \leqslant \rho_{2}(t) \quad \forall \, r \in[t-1, t],\, \tau \leqslant \tau_{1}(t, \widehat{D}),\, u_{\tau} \in D(\tau).
\label{4.18-2}
\end{equation}

From the proof of Theorem \ref{th3.1-2}, we can conclude that ${u_k} \to u(t;\tau ,{u_\tau })$ weakly in ${L^\infty }( {\tau ,T;{{\cal H}_t}(\Omega )} ) \cap {L^2}( {\tau ,T;H_0^1(\Omega )} )$ for all $r \in [t - 1,t]$. Then according to (\ref{4.13-2}) and (\ref{4.18-2}), the inequality (\ref{4.14-2}) follows directly. $\hfill$$\Box$

Now we will use the energy method to prove that the process $\left\{ {U(t,\tau ){\} _{t \ge \tau }}} \right.$ is pullback $\mathcal{D}_{\sigma}^{\mathcal H_{t}}$-asymptotically compact.
\begin{Lemma}
Under the assumptions of Lemma $\ref{lem4.2-2}$, the process $\left\{ {U(t,\tau ){\} _{t \ge \tau }}} \right.$ on $\mathcal H_{t}(\Omega)$ is pullback $\mathcal{D}_{\sigma}^{\mathcal H_{t}}$-asymptotically compact.
\label{lem4.4-2}
\end{Lemma}
$\mathbf{Proof.}$ Let $t \in \mathbb{R}$, $\widehat{D} \in \mathcal{D}_{\sigma}^{\mathcal H_{t}}$, $\left\{\tau_{k}\right\} \subset(-\infty, t-2]$ with $\tau_{k} \rightarrow-\infty$, and $u_{\tau_{k}} \in D\left(\tau_{k}\right)$ for all $k \in \mathbb N^{+}$. According to Definition \ref{def2.4-2}, if it is proved that the sequence ${\left\{u\left(t ; \tau_{k}, u_{\tau_{k}}\right)\right\}}$ is relatively compact in $\mathcal H_{t}$, then $u^{k}(t)={\left\{u\left(t ; \tau_{k}, u_{\tau_{k}}\right)\right\}}$ is pullback $\mathcal{D}_{\sigma}^{\mathcal H_{t}}$-asymptotically compact.

By similar calculations in the proof of Theorems \ref{th3.1-2} and \ref{th3.2-2} and Lemma \ref{lem4.3-2}, it follows that there exists $\tau_{1}(t,\widehat{D})<t-2$ such that ${\{ {u^k}\} _{k \ge {k_1} \ge 1}}$ is bounded in ${L^\infty }( {t-2 ,t;{{\cal H}_t}(\Omega )} ) \cap {L^2}( {t-2 ,t;H_0^1(\Omega )} )$, ${\{ f({u^k})\} _{k \ge {k_1} \ge 1}}$ is bounded in $L^{q}(t-2, t ; L^{q}(\Omega))$ and ${\{ {\partial _t}{u^k}\} _{k \ge {k_1} \ge 1}}$ is bounded in $L^{2}(t-2, t ; {{\cal H}_t}(\Omega ))$. Then, by the Aubin-Lions compactness lemma, we derive that there exists $u \in L^{2}(t-2, t ; {{\cal H}_t}(\Omega ))$ with ${\partial _t}{u} \in L^{2}(t-2, t ; {{\cal H}_t}(\Omega ))$ such that
\begin{equation}
\quad \quad u^{k} \rightharpoonup u \quad \text { weakly-star in } L^{\infty}(t-2, t ;\mathcal H_{t}(\Omega));
\label{4.19-2}
\end{equation}
\begin{equation}
u^{k} \rightharpoonup u \quad \text { weakly in } L^{2}(t-2, t ;H_0^1(\Omega ));
\label{4.20-2}
\end{equation}
\begin{equation}
\quad f(u^{k}) \rightharpoonup f(u) \quad \text { weakly in } L^{q}(t-2, t ; L^{q}(\Omega));
\label{4.21-2}
\end{equation}
\begin{equation}
\quad\quad\quad\quad a(l(u^{k})) u^{k} \rightharpoonup a(l(u)) u \quad \text { weakly in } L^{2}(t-2, t ; H_{0}^{1}(\Omega));
\label{4.22-2}
\end{equation}
\begin{equation}
\quad {\partial _t}{u^k} \rightharpoonup {\partial _t}u \quad \text { weakly in } L^{2}(t-2, t ;\mathcal H_{t}(\Omega));
\label{4.23-2}
\end{equation}
\begin{equation}
u^{k} \rightarrow u \quad a. e.\quad(x, t) \in \Omega \times[\tau,+\infty),
\label{4.24-2}
\end{equation}
for any $\tau_{k} \leq \tau_{1}(t,\widehat{D})$. Then by the similar proof to Theorem \ref{th3.1-2}, we obtain that $u \in C([t-2,t]; \mathcal H_{t}(\Omega))$, and from $(\ref{4.19-2})-(\ref{4.24-2})$, we can derive that $u$ satisfies Definition \ref{def3.1-2} in $(t-2,t)$.

In order to prove the lemma, we only need to check that ${u^k} \to u$ strongly in $C([t-2,t]; \mathcal H_{t}(\Omega))$. We establish it by contradiction. Suppose there is a constant $\mu  > 0$ and a sequence $t_{k} \in [t-2,t]$ such that when ${t_k} \to {t_*}$,
\begin{equation}
|u^{k}(t_{k})-u(t_{*})| \geq \mu, \quad \forall \, k \in {\mathbb N^{+}}.
\label{4.25-2}
\end{equation}

By similar calculations to (\ref{3.11-2}), we derive
\begin{equation}
\|e(s)\|_{\mathcal{H}_{t}}^{2} \leqslant\|e(r)\|_{\mathcal{H}_{t}}^{2}+2 C_{0}(s-r)+\frac{2 \xi^{2}}{m \lambda_{1}} \int_{r}^{s}\|h(x, z)\|^{2} d z, \quad \forall \, t-2 \leqslant r \leqslant s \leqslant t,
\label{4.26-2}
\end{equation}
where the function $e$ can be replaced by $u$ or $u^{k}$.

Now we define the following functions
\begin{equation}
Q_{k}(s)=\|u^{k}(s)\|_{\mathcal H_{t}}^{2}-2 C_{0}s -\frac{2 \xi^{2}}{m \lambda_{1}} \int_{t-2}^{s}\|h(x, r)\|^{2} d r
\label{4.27-2}
\end{equation}
and
\begin{equation}
Q(s)=\|u(s)\|_{\mathcal H_{t}}^{2}-2 C_{0}s -\frac{2 \xi^{2}}{m \lambda_{1}} \int_{t-2}^{s}\|h(x, r)\|^{2} d r.
\label{4.28-2}
\end{equation}

From the smoothness of $u$ and $u^{k}$, we obtain that the functions $Q_{k}(s)$ and $Q(s)$ are continuous on $[t-2, t]$. Then using the above inequality, we obtain that $Q_{k}(s)$ and $Q(s)$ are non-increasing on $[t-2, t]$ and from $(\ref{4.19-2})-(\ref{4.24-2})$, we can derive
\begin{equation}
Q_{k}(s) \rightarrow Q(s) \quad \text { a.e. } s \in(t-2, t).
\label{4.29-2}
\end{equation}

Assume that there is a sequence $\left\{\tilde{t}_{n}\right\}\left(t-2, t_{*}\right)$ such that $\tilde{t}_{n} \rightarrow t_{*}$ when $k \rightarrow \infty$ and the above convergence holds. Since the function $Q(s)$ is continuous on $[t-2, t]$, there exists a positive integer $n(\tilde{\varepsilon}) \geqslant 1$ such that
\begin{equation}
\left|Q\left(\tilde{t}_{n}\right)-Q\left(t_{*}\right)\right|<\frac{\tilde{\varepsilon}}{2}, \quad \forall \, n \geqslant n(\tilde{\varepsilon}),
\label{4.30-2}
\end{equation}
where $\tilde \varepsilon$ is a positive constant.

From (\ref{4.29-2}), we can derive that there exists a constant $k(\tilde{\varepsilon}) \in \mathbb{N}^{+}$ such that $t_{k} \geqslant \tilde{t}_{n(\tilde{\varepsilon})}$ and $\left|Q_{k}\left(\tilde{t}_{n(\tilde{\varepsilon})}\right)-Q\left(\tilde{t}_{n(\tilde{\varepsilon})}\right)\right|<\frac{\tilde{\varepsilon}}{2}$ for all $k \geqslant k(\tilde{\varepsilon})$. Furthermore, notice that $Q_{k}(s)$ is a non-increasing function, then by the Cauchy inequalities, we can conclude that for all $k \geqslant k(\tilde{\varepsilon})$ the following inequality holds
\begin{equation}
\begin{aligned}
Q_{k}\left(t_{k}\right)-Q\left(t_{*}\right) & \leqslant |Q_{k}(\tilde{t}_{n(\tilde{\varepsilon})})-Q(t_{*})|\\
& \leqslant |Q_{k} (\tilde{t}_{n(\tilde{\varepsilon})})-Q (\tilde{t}_{n(\tilde{\varepsilon})}) |+ |Q (\tilde{t}_{n(\tilde{\varepsilon})})-Q(t_{*})| \\
& \leqslant \frac{\tilde{\varepsilon}}{2}+\frac{\tilde{\varepsilon}}{2}=\tilde{\varepsilon}.
\end{aligned}
\label{4.31-2}
\end{equation}

Since the arbitrariness of $\tilde \varepsilon>0 $, it follows that $\lim\limits_{k \rightarrow \infty} \sup Q_{k}\left(t_{k}\right) \leqslant Q\left(t_{*}\right)$ as $\tilde \varepsilon \rightarrow 0.$ Then we can conclude that $\lim\limits _{k \rightarrow \infty} \sup \|u^{k}(t_{k})\|_{\mathcal H_{t}} \leqslant \|u(t_{*}) \|_{\mathcal H_{t}}$, which contradicts with $(\ref{4.25-2})$. Therefore, the proof is complete. $\hfill$$\Box$

From the above proofs, it follows the following theorem about the existence of the minimal time-dependent pullback $\mathcal{D}_{\sigma}^{\mathcal H_{t}}$-attractors.

\begin{Theorem}\label{th4.1-2}
Under the assumptions of Theorems $\ref{th3.1-2}$ and $\ref{th3.2-2}$ and assume that the function $h(x, t)$ satisfies $(\ref{4.11-2})$, then there exists the minimal time-dependent pullback $\mathcal D_{F}^{\mathcal H_{t}}$-attractor $\mathcal A_{\mathcal D_{F}^{\mathcal H_{t}}}=\{A_{\mathcal D_{F}^{\mathcal H_{t}}}( t):t \in \mathbb{R}\}$ and the minimal time-dependent pullback $\mathcal D_{\sigma}^{\mathcal H_{t}}$-attractors $\mathcal A_{\mathcal D_{\sigma}^{\mathcal H_{t}}}=\{A_{\mathcal D_{\sigma}^{\mathcal H_{t}}}( t):t \in \mathbb{R}\}$ for the process $\{U(t, \tau)\}_{t \geqslant \tau}$ on $\mathcal H_{t}(\Omega)$  of problem $(\ref{1.1-2})$. Moreover, the family $\mathcal A_{\mathcal D_{\sigma}^{\mathcal H_{t}}}$ belongs to $\mathcal D_{\sigma}^{\mathcal{H}_{t}}$ and for any $t \in \mathbb{R}$ the following relationships hold:
\begin{equation}
A_{\mathcal{D}_{F}^{\mathcal{H}_{t}}}(t) \subset A_{\mathcal{D}_{\sigma}^{\mathcal{H}_{t}}}(t) \subset \bar{B}_{\mathcal{H}_{t}}\left(0, \rho_{\xi}^{\frac{1}{2}}(t)\right).
\label{4.32-2}
\end{equation}
In addition, if $h(x, t)$ satisfies
$$
\sup _{r \leq 0}\left(e^{-\sigma r} \int_{-\infty}^{r} e^{\sigma s}\|h(x, s)\|^{2} d s\right)<+\infty,
$$
then $A_{\mathcal {D}_{F}^{\mathcal H_{t}}}(t)= A_{\mathcal {D}_{\sigma}^{\mathcal H_{t}}}(t)$ for any $t \in \mathbb{R}$.
\end{Theorem}
$\mathbf{Proof.}$ From Definition \ref{def2.9-2}, Theorems $\ref{th3.1-2}-\ref{th3.2-2}$ and Lemmas $\ref{lem4.1-2}-\ref{lem4.4-2}$, it follows the existence of above minimal time-dependent pullback $\mathcal{D}_{\sigma}^{\mathcal H_{t}}$-attractors. Besides, from Lemma \ref{th2.1-2} and Corollaries $\ref{re2.3-2}-\ref{re2.4-2}$, we can derive that (\ref{4.32-2}) and $A_{\mathcal {D}_{F}^{\mathcal H_{t}}}(t)= A_{\mathcal {D}_{\sigma}^{\mathcal H_{t}}}(t)$. $\hfill$$\Box$

\section{Existence of time-dependent pullback attractors $\left\{A_{\xi}(t)\right\}_{t \in \mathbb R}$ and upper semicontinuity of
$\left\{A_{\xi}(t)\right\}_{t \in \mathbb R}$ and the global attractor $A$}
\ \ \ \ In this section, we will establish the existence of time-dependent pullback attractors $\left\{A_{\xi}(t)\right\}_{t \in \mathbb R}$ and the upper semicontinuity of $\left\{A_{\xi}(t)\right\}_{t \in \mathbb R}$ and the global attractor $A$ of equation (\ref{1.1-2}) with $\xi = 0$.

First of all, we decompose the solution $U_{\xi}(t, \tau) u_{\tau}=u(t)$ of problem (\ref{1.1-2}) with initial data $u_{\tau} \in \mathcal{H}_{t}(\Omega)$ as follows:
$$
U_{\xi}(t, \tau) u_{\tau}=U_{1, \xi}(t, \tau) u_{\tau}+U_{2, \xi}(t, \tau) u_{\tau},
$$
where $U_{1, \xi}(t, \tau) u_{\tau}=v(t)$ and $U_{2, \xi}(t, \tau) u_{\tau}=g(t)$ solve, respectively,
\begin{equation}
\left\{\begin{array}{ll}
v_{t}-\varepsilon (t) \Delta v_{t}-a(l(u)) \Delta v=f_{0}(v) & \text { in } \Omega\times(\tau, \infty),  \\
v=0 & \text { on } \partial \Omega \times(\tau, \infty), \\
v(x, \tau)=u_{\tau}(x) &\,\,  x \in \Omega,
\end{array}\right.\label{5.1-2}
\end{equation}
and
\begin{equation}
\left\{\begin{array}{ll}
g_{t}-\varepsilon(t) \Delta g_{t}-a(l(u)) \Delta g=f(u)-f_{0}(v)+\xi h(x, t) & \text { in } \Omega \times(\tau, \infty),   \\
g=0 & \text { on } \partial \Omega \times(\tau, \infty),  \\
g(x, \tau)=0 &\,\, x \in \Omega .
\end{array}\right.\label{5.2-2}
\end{equation}
\begin{Lemma}
Under the assumptions of Theorems $\ref{th3.1-2}$ and $\ref{th3.2-2}$ and Lemma $\ref{lem4.1-2}$, if for any $t \in \mathbb R$ the function $h(x,t)$ satisfies $(\ref{4.11-2})$, then for any bounded set $B \subset \mathcal{H}_{t}(\Omega)$, there exists $T(t, B)>0$ such that
\begin{equation}
\left\|U_{\xi}(t, t-\tau) u_{t-\tau}\right\|_{\mathcal{H}_{t}}^{2} \leq \rho_{\xi}(t) \quad \text { for all } \tau \geq T(t, B) \text { and all } u_{t-\tau} \in B,
\label{5.3-2}
\end{equation}
where $\rho_{\xi}(t)$ is the same as $(\ref{4.12-2})$.
\label{lem5.1-2}
\end{Lemma}
$\mathbf{Proof.}$  From Lemmas $\ref{lem4.1-2}$ and $\ref{lem4.2-2}$, the desired result follows easily. $\hfill$$\Box$

Let ${D_\xi }(t) = \{ u \in {{\cal H}_t}(\Omega )\left| \|u\|_{\mathcal{H}_{t}}^{2} \leq \rho_{\xi}(t) \right.\}$. It is easy to check that the family $\mathcal{D}_{\xi}=\left\{D_{\xi}(t)\right\}_{t \in \mathbb{R}}$ is pullback absorbing in $\mathcal{H}_{t}(\Omega)$. Moreover, we derive
$$
\lim _{t \rightarrow-\infty} e^{\sigma t} \rho_{\xi}(t)=0 \quad \text { for any } \xi>0 .
$$
\begin{Lemma}
 Let $\rho_{\xi}(t), D_{\xi}(t)$ given as the above. For any $t \in \mathbb R$, the solution $v(t)$ of problem $(\ref{5.1-2})$ satisfies
\begin{equation}
\|v(t)\|_{\mathcal{H}_{t}}^{2}=\left\|U_{1, \xi}(t, t-\tau) u_{t-\tau}\right\|_{\mathcal{H}_{t}}^{2} \leq\left(1+2 \lambda_{1}\right) e^{-\sigma \tau} \rho_{\xi}(t-\tau)
\label{5.4-2}
\end{equation}
for all $\tau \geq 0$, $u_{t-\tau} \in D_{\xi}(t-\tau)$ and $0<\sigma<\min \left\{\eta, \frac{-\varepsilon^{\prime}(t)}{\varepsilon(t)}, \frac{1}{\left(1+\lambda_{1}\right) \lambda_{1}^{-1}}\right\}$ with $0<\eta<m\lambda_{1}$.
\label{lem5.2-2}
\end{Lemma}
$\mathbf{Proof.}$ Multiplying $(\ref{5.1-2})_{1}$ by $v$, then we arrive at
\begin{equation}
\frac{d}{d t}(\|v\|^{2}+\varepsilon(t)\|\nabla v\|^{2})+(2 a(l(u))-\varepsilon^{\prime}(t))\|\nabla v\|^{2}=2(f_{0}(v), v).
\label{5.5-2}
\end{equation}

Then from $(\ref{1.9-2})$, we can derive
\begin{equation}
\frac{d}{d t}\left(\|v\|^{2}+\varepsilon(t)\|\nabla v\|^{2}\right)+\frac{\lambda_{1}}{1+\lambda_{1}}\left(\|v\|^{2}+\varepsilon(t)\|\nabla v\|^{2}\right) \leqslant 0
\label{5.6-2}
\end{equation}
which, by taking $0<\sigma<\min \left\{\eta, \frac{-\varepsilon^{\prime}(t)}{\varepsilon(t)}, \frac{1}{\left(1+\lambda_{1}\right) \lambda_{1}^{-1}}\right\}$ with $0<\eta<m\lambda_{1}$, gives
\begin{equation}
\frac{d}{d t}(e^{\sigma t}(\|v\|^{2}+\varepsilon (t)\|\nabla v \|^{2}))+\left(\frac{1}{\left(1+\lambda_{1}\right) \lambda_{1}^{-1}}-\sigma\right) e^{\sigma t}\left(\|v\|^{2}+\varepsilon(t)\|\nabla v\|^{2}\right) \leqslant 0.
\label{5.7-2}
\end{equation}

Integrating (\ref{5.7-2}) over $\left[t_{0}-\tau, t_{0}\right]$, then by Lemmas \ref{lem4.1-2} and \ref{lem5.1-2}, we conclude
\begin{equation}
\begin{aligned}
\left\|v\left(t_{0}\right)\right\|^{2}+\varepsilon(t)\left\|\nabla v\left(t_{0}\right)\right\|^{2} & \leqslant e^{-\sigma \tau}(\|v_{t_{0}-\tau}\|^{2}+\varepsilon(t_{0}-\tau)\|\nabla v_{t_{0}-\tau}\|^{2}) \\
& \leqslant e^{-\sigma \tau}(\|u_{t_{0}-\tau} \|^{2}+\varepsilon(t_{0}-\tau)\| \nabla u_{t_{0}-\tau}\|^{2}).\\
& \leqslant(1+2 \lambda_{1}) e^{-\sigma \tau} \rho_{\xi}(t-\tau),
\end{aligned}
\label{5.8-2}
\end{equation}
which completes the proof. $\hfill$$\Box$
\begin{Lemma}\label{lem5.3-2}
For any $t \in \mathbb{R}$, there exists $T(t, \mathcal{D}_{\xi})>0$ and $R_{\xi}(t)>0$, such that
\begin{equation}
{\|U_{2, \xi}(t, t-\tau) u_{t-\tau}\|}^{2}_{\mathcal{H}_{t}^{1+\alpha}}=\|A^{\frac{\alpha}{2}} g\|^{2}+\varepsilon(t)\|A^{\frac{1+\alpha}{2}} g\|^{2} \leqslant R_{\xi}(t)
\label{5.9-2}
\end{equation}
for any $\tau  > T(t,\mathcal D_{\xi})$ and $u_{t-\tau} \in \mathcal D_{\xi}(t-\tau)$, where $R_{\xi}(t)=C \rho_{\xi}^{2 p+2}(t)$ and $0<\alpha<\min \{1, \frac{4-(N-2) p}{2}\}$.
\end{Lemma}
$\mathbf{Proof.}$ Multiplying $(\ref{5.2-2})_{1}$ by $e^{\delta t} A^{\alpha} g$, then we obtain
\begin{equation}
(g_{t}, e^{\delta t} A^{\alpha} g)-(\varepsilon(t) \Delta g_{t}, e^{\delta t} A^{\alpha} g)-(a(l(u)) \Delta g, e^{\delta t} A^{\alpha} g)=(f(u)-f_{0}(v)+\xi h(x,t), e^{\delta t} A^{\alpha} g),
\label{5.10-2}
\end{equation}
where $0<\delta<\min \left\{(p+1) \sigma, \frac{16 m-2-\xi}{8\left(\lambda_{1}^{-1}+\varepsilon(t)\right)}\right\}$.

By simple calculations, we derive the following equalities
\begin{equation}
\begin{aligned}
(g_{t}, e^{\delta t} A^{\alpha} g) &=\frac{1}{2} \frac{d}{d t}(e^{\delta t}\|A^{\frac{\alpha}{2}} g\|^{2})-\frac{1}{2} \delta e^{\delta t} \| A^{\frac{\alpha}{2}} g \|^{2},
\end{aligned}
\label{5.11-2}
\end{equation}
\begin{equation}
\begin{aligned}
-(\varepsilon (t) \Delta g_{t}, e^{\delta t} A^{\alpha} g) =\frac{1}{2} \frac{d}{d t}(\varepsilon(t) e^{\delta t}\|A^{\frac{1+\alpha}{2}} g\|^{2})-\frac{1}{2} \varepsilon^{\prime}(t) e^{\delta t}\|A^{\frac{1+\alpha}{2}} g\|^{2}-\frac{1}{2} \delta \varepsilon(t) e^{\delta t}\|A^{\frac{1+\alpha}{2}} g\|^{2}
\end{aligned}
\label{5.12-2}
\end{equation}
and
\begin{equation}
-(a(l(u)) \Delta g, e^{\delta t} A^{\alpha} g)=a (l(u)) e^{\delta t}\|A^{\frac{1+\alpha}{2}}g\|^{2}.
\label{5.13-2}
\end{equation}

From (\ref{1.7-2}), we conclude
\begin{equation}
(f(u)-f_{0}(v), e^{\delta t} A^{\alpha} g)=e^{\delta t}(f(u)-f(v), A^{\alpha} g)+e^{\delta t}(f_{1}(v), A^{\alpha} g).
\label{5.14-2}
\end{equation}

Inserting $(\ref{5.11-2})-(\ref{5.14-2})$ into (\ref{5.10-2}), we obtain
\begin{equation}
\begin{aligned}
&\frac{d}{d t}(e^{\delta t}(\|A^{\frac{\alpha}{2}} g\|^{2}+\varepsilon(t)\|A^{\frac{1+\alpha}{2}} g\|^{2}))-\varepsilon^{\prime}(t) e^{\delta t}\|A^{\frac{1+\alpha}{2}} g\|^{2}+2 a(l(u)) e^{\delta t}\|A^{\frac{1+\alpha}{2}} g\|^{2} \\
&=\delta e^{\delta t}(\|A^{\frac{\alpha}{2}} g\|^{2}+\varepsilon(t)\|A^{\frac{1+\alpha}{2}} g\|^{2})+2 e^{\delta t}(f(u)-f(v), A^{\alpha} g)+2 e^{\delta t}(f_{1}(v), A^{\alpha} g)\\
& +2 \xi e^{\delta t}(h, A^{\alpha} g).
\end{aligned}
\label{5.15-2}
\end{equation}

By similar calculations to  Lemma 3.3 in \cite{wlq.2}, noticing that $p \le \frac{4}{{N - 2}}$, from assumptions (\ref{1.7-2}),  (\ref{1.10-2}) and Lemma \ref{1m2.1-2}, then assuming $0<\alpha<\min \{1, \frac{4-(N-2) p}{2}\}$, we can derive the following inequalities
\begin{equation}
|e^{\delta t}(f(u)-f(v), A^{\alpha} g)| \leqslant C e^{\delta t}(\|A^{\frac{1}{2}} u\|^{2 p}+\|A^{\frac{1}{2}} v\|^{2 p}+1)\|A^{\frac{1}{2}} g\|^{2})+\frac{1}{16} e^{\delta t}\|A^{\frac{1+ \alpha}{2}} g\|^{2},
\label{5.16-2}
\end{equation}
\begin{equation}
|e^{\delta t}(f_{1}(v), A^{\alpha} g)| \leqslant Ce^{\delta t}(1+\|A^{\frac{1}{2}} v\|^{2p+2})+\frac{1}{16} e^{\delta t}\|A^{\frac{1 + \alpha}{2}} g\|^{2}
\label{5.17-2}
\end{equation}
and
\begin{equation}
|(\xi h(x, t), A^{\alpha} g)| \leq 4 \xi\|h(x, t)\|^{2}+\frac{1}{16}\xi\|A^{\frac{1+\alpha}{2}} g\|^{2}.
\label{5.18-2}
\end{equation}

Inserting $(\ref{5.16-2})-(\ref{5.18-2})$ into (\ref{5.15-2}), we derive
\begin{equation}
\begin{aligned}
&\frac{d}{d t}(e^{\delta t}(\|A^{\frac{\alpha}{2}} g\|^{2}+\varepsilon(t)\|A^{\frac{1+ \alpha}{2}} g\|^{2}))-\varepsilon^{\prime}(t) e^{\delta t}\|A^{\frac{1+ \alpha}{2}} g\|^{2}.\\
&\leqslant Ce^{\delta t}(\|A^{\frac{1}{2}} u\|^{2}+\|A^{\frac{1}{2}}v \|^{2}+1)\|A^{\frac{1}{2}}g\|^{2}+C e^{\delta t}(1+\|A^{\frac{1}{2}} v\|^{2 p+2})+c\xi e^{\delta t}\|h(x, t)\|^{2}.
\end{aligned}
\label{5.19-2}
\end{equation}

From $u=v+g$, we conclude
\begin{equation}
\begin{aligned}
&(\|A^{\frac{1}{2}} u\|^{2 p}+\|A^{\frac{1}{2}} v\|^{2 p}+1)\|A^{\frac{1}{2}} g\|^{2} \\
&\leqslant C(\|A^{\frac{1}{2}} u\|^{2 p}\|A^{\frac{1}{2}} u-A^{\frac{1}{2}} v\|^{2}+\|A^{\frac{1}{2}} v\|^{2 p}\|A^{\frac{1}{2}} u-A^{\frac{1}{2}} v\|^{2}+1) \\
&\leqslant C(\|A^{\frac{1}{2}} u\|^{2 p+2}+\|A^{\frac{1}{2}} v\|^{2 p+2}+1).
\end{aligned}
\label{5.20-2}
\end{equation}

Inserting (\ref{5.20-2}) into (\ref{5.19-2}) and by the Poincar$\acute{e}$ inequality, we conclude
\begin{equation}
\begin{aligned}
&\frac{d}{d t}(e^{\delta t}(\|A^{\frac{\alpha}{2}} g\|^{2}+\varepsilon(t)\|A^{\frac{1+\alpha}{2}} g\|^{2}))-\varepsilon^{\prime}(t) e^{\delta t}\|A^{\frac{1+\alpha}{2}} g\|^{2} \\
&\leqslant C e^{\delta t}(\|A^{\frac{1}{2}} u\|^{2 p+2}+\|A^{\frac{1}{2}} v\|^{2 p+2}+1)+C \xi e^{\delta t}\|h(x, t)\|^{2}.
\end{aligned}
\label{5.21-2}
\end{equation}

From (\ref{4.1-2}), it follows that
\begin{equation}
\|u\|_{\mathcal{H}_{t}}^{2} \leqslant C e^{-\sigma \tau}\left\|u_{t-\tau}\right\|_{\mathcal{H}_{t}}^{2}+C\left(\xi e^{-\sigma t} \int_{-\infty}^{t} e^{\sigma s}\|h(x, s)\|^{2} d s+1\right).
\label{5.22-2}
\end{equation}

Multiplying (\ref{5.22-2}) by $e^{\sigma t}$, then we obtain
\begin{equation}
e^{\sigma t}\|u\|_{\mathcal H_{t}}^{2} \leqslant C e^{\sigma (t_{0}-\tau)}\|u_{t_{0}-\tau}\|_{\mathcal H_{t}}^{2}+C \xi \int_{-\infty}^{t_{0}} e^{\sigma t}\| h(x, t) \|^{2} d t+C e^{\sigma t}.
\label{5.23-2}
\end{equation}

Calculating the $p+1$ power on both sides of (\ref{5.23-2}),  we conclude
\begin{equation}
\begin{aligned}
e^{(p+1) \sigma t}\|u\|_{\mathcal{H}_{t}}^{2 p+2} & \leqslant C e^{(p+1) \sigma}\left\|u_{t_{0}-\tau}\right\|_{\mathcal{H}_{t}}^{2 p+2}+C \xi\left(\int_{-\infty}^{t_{0}} e^{\sigma t}\|h(x, t)\|^{2} d t\right)^{2 p+2} \\
&+C e^{(p+1) \sigma t}+C .
\end{aligned}
\label{5.24-2}
\end{equation}

Multiplying (\ref{5.24-2}) by $ e^{(\delta-(p+1) \sigma) t}$, then we derive
\begin{equation}
\begin{aligned}
e^{\delta t}\|u\|_{\mathcal{H}_{t}}^{2 p+2} &\leqslant  C e^{(p+1) \sigma\left(t_{0}-\tau\right)} e^{(\delta-(p+1) \sigma) t}\left\|u_{t_{0}-\tau}\right\|_{\mathcal{H}_{t}}^{2 p+2}\\
&+C \xi e^{(\delta-(p+1) \sigma) t}\left(\int_{-\infty}^{t_{0}} e^{\sigma t}\|h(x, t)\|^{2} d t\right)^{2 p+2}+C e^{\delta t}+C e^{(\delta-(p+1) \sigma) t}.
\end{aligned}
\label{5.25-2}
\end{equation}

Integrating (\ref{5.25-2}) over $\left[t_{0}-\tau, t_{0}\right]$, we conclude
\begin{equation}
\begin{aligned}
&\int_{-\infty}^{t_{0}} e^{\delta t}{\| u(t) \|_{\mathcal{H}_{t}}^{2 p+2}} d t \leqslant C e^{(p+1) \sigma\left(t_{0}- \tau\right)} e^{(\delta-(p+1)\sigma) t} \| u_{t_{0}-\tau}{ \|_{\mathcal{H}_{t}}^{2 p+2}}\\
& \,\, \quad \quad \quad \quad \quad \quad \quad \quad +C \xi e^{(\delta-(p+1)\sigma )t_{0} }\left(\int_{-\infty}^{t_{0}} e^{\sigma t}\|h(x, t)\|^{2} d t\right)^{p+1}+C e^{\delta t_{0}}+C e^{(\delta-(p+1) \sigma) t_{0}}.
\end{aligned}
\label{5.26-2}
\end{equation}

By(\ref{5.8-2}) and the Poincar$\acute{e}$ inequality, we conclude
\begin{equation}
e^{\sigma t}\|v(t)\|_{\mathcal{H}_{t}}^{2} \leqslant C e^{\sigma\left(t_{0}-\tau\right)}\left\|u_{t_{0}-\tau }\right\|_{\mathcal{H}_{t}}^{2}.
\label{5.27-2}
\end{equation}

Then calculating the $p+1$ power on both sides of (\ref{5.27-2}), we obtain
\begin{equation}
e^{(p+1) \sigma t}\|v(t)\|_{\mathcal H_{t}}^{2} \leqslant Ce^{(p+1) \sigma(t_{0}-\tau)} \|u_{t_{0}-\tau} \|_{\mathcal H_{t}}^{2}.
\label{5.28-2}
\end{equation}

Calculating similarly to deriving (\ref{5.25-2}) from (\ref{5.24-2}), we can derive
\begin{equation}
e^{\delta t}\|v(t)\|_{\mathcal H_{t}}^{2 p+2} \leqslant C e^{(\delta-(p+1) \sigma) t} e^{(p+1) \sigma\left(t_{0}-\tau\right)} \| u_{t_{0}-\tau}\|_{\mathcal H_{t}} ^{2 p+2}.
\label{5.29-2}
\end{equation}

Then integrating (\ref{5.29-2}) over $\left[t_{0}-\tau, t_{0}\right]$, we obtain
\begin{equation}
\int_{-\infty}^{t_{0}} e^{\sigma t}\|v(t)\|_{\mathcal H_{t}}^{2 p+2} d t \leqslant C e^{\sigma  t_{0}} e^{-(p+1) \sigma \tau}\left\|u_{t_{0}-\tau}\right\|_{\mathcal H_{t}}^{2 p+2}.
\label{5.30-2}
\end{equation}

Integrating (\ref{5.21-2}) over $\left[t_{0}-\tau, t_{0}\right]$ and noting that $\varepsilon(t)$ is a decreasing function, we obtain
\begin{equation}
\begin{aligned}
e^{\delta t_{0}}(\|A^{\frac{\alpha}{2}} g(t_{0})\|^{2}+\varepsilon(t_{0})\|A^{\frac{1+\alpha}{2}} g(t_{0})\|^{2}) & \leqslant e^{\delta(t_{0}-\tau)}(\|A^{\frac{\alpha}{2}} g_{t_{0}-\tau}\|^{2}+\varepsilon(t_{0}-\tau)\|A^{\frac{1+\alpha}{2}} g_{t_{0}-\tau}\|^{2}) \\
&+C \int_{t_{0}-\tau}^{t_{0}} e^{\delta t}(\|A^{\frac{1}{2}}u\|^{2 p+2}+\| A^{\frac{1}{2}}v\|^{2 p+2}+1) d t \\
&+C\xi \int_{t_{0}-\tau}^{t_{0}} e^{\delta t}\|h(x, t)\|^{2} d t+C.
\end{aligned}
\label{5.31-2}
\end{equation}

Then we can derive
\begin{equation}
\begin{aligned}
&\|A^{\frac{\alpha}{2}} g(t_{0})\|^{2}+\varepsilon(t_{0})\|A^{\frac{1+\alpha}{2}} g(t_{0})\|^{2} \leqslant C e^{-\delta t_{0}} \int_{t_{0}-\tau}^{t_{0}}e^{\delta t}(\|A^{\frac{1}{2}}u\|^{2 p+2}+\| A^{\frac{1}{2}}v\|^{2 p+2}+1) d t\\
&\qquad \qquad\qquad\qquad\qquad\qquad\qquad+C \xi e^{-\delta t_{0}} \int_{t_{0}-\tau}^{t_{0}} e^{\delta t}\|h(x, t)\|^{2} d t+C.
\end{aligned}
\label{5.32-2}
\end{equation}

Multiplying (\ref{5.26-2}) by $ e^{-\delta t_{0}}$ and by $\delta >(p+1)\sigma$, we obtain
\begin{equation}
\begin{aligned}
e^{-\delta t_{0}} \int_{-\infty}^{t_{0}}& e^{\delta t}\|u(t)\|_{\mathcal {H}_{t}}^{2 p+2} d t \leq C e^{-\delta t_{0}} e^{(p+1) \sigma\left(t_{0}-\tau\right)} e^{(\delta-(p+1) \sigma) t}\left\|u_{t_{0}-\tau}\right\|_{\mathcal {H}_{t}}^{2 p+2} \\
&+C \xi e^{(\delta-(p+1) \sigma) t_{0}} e^{-\delta t_{0}}\left(\int_{-\infty}^{t_{0}} e^{\sigma t}\|h(x, t)\|^{2} d t\right)^{p+1}+C e^{-(p+1) \sigma t_{0}}+C \\
& \leq C e^{-(p+1) \sigma \tau}\left\|u_{t_{0}-\tau}\right\|_{\mathcal {H}_{t}}^{2 p+2}+C \xi e^{-(p+1) \sigma t_{0}}\left(\int_{-\infty}^{t_{0}} e^{\sigma t}\|h(x, t)\|^{2} d t\right)^{p+1}+C \\
& \leqslant C e^{-\delta \tau }\|u_{t_{0}-\tau} \|_{\mathcal H_{t}}^{2 p+2}+C \xi e^{-\delta t_{0}}\left(\int_{-\infty}^{t_{0}} e^{\sigma t}\|h(x, t)\|^{2} d t\right)^{p+1}+C.
\end{aligned}
\label{5.33-2}
\end{equation}

Besides, multiplying (\ref{5.30-2}) by $ e^{-\delta t_{0}}$, we conclude
\begin{equation}
e^{-\delta t_{0}} \int_{-\infty}^{t_{0}} e^{\sigma t}\|v(t)\|_{\mathcal H_{t}}^{2 p+2} d t \leqslant e^{-(p+1) \sigma \tau}\left\|u_{t_{0}-\tau}\right\|_{\mathcal H_{t}}^{2 p +2} \leqslant e^{-\delta \tau} \| u_{t_{0}-\tau }\|{_{\mathcal H_{t}}^{2 p+2}}.
\label{5.34-2}
\end{equation}

Then inserting (\ref{5.33-2}) and (\ref{5.34-2}) into (\ref{5.32-2}) and by Lemma \ref{1m2.1-2}, we derive
\begin{equation}
\begin{aligned}
&\|A^{\frac{\alpha}{2}} g(t_{0})\|^{2}+\varepsilon(t_{0})\|A^{\frac{1+\alpha}{2}} g(t_{0})\|^{2} \leqslant C+C e^{-\delta t_{0}} \int_{-\infty}^{t_{0}} e^{\delta t}\|A^{\frac{1}{2}}u\|^{2 p+2} d t\\
&+C e^{-\delta t_{0}} \int_{-\infty}^{t_{0}} e^{\delta t}\|A^{\frac{1}{2}} v\|^{2 p+2} d t+C e^{-\delta t_{0}}{\int_{t_{0}-\tau}^{t_{0}} e^{\delta t} d t}+C\xi e^{-\delta t_{0}} \int_{-\infty}^{t_{0}} e^{\delta t}\|h(x, t)\|^{2} d t\\
&\leq C e^{-\delta \tau}\left\|u_{t_{0}-\tau}\right\|_{\mathcal{H}_{t}}^{2 p+2}+C \xi e^{-\delta t_{0}} \int_{-\infty}^{t_{0}} e^{\delta t}\|h(x, t)\|^{2} d t+C .
\end{aligned}
\label{5.35-2}
\end{equation}

From $p \le \frac{4}{{N - 2}}$, it follows that $p+1 \le 5$. Then by Lemma \ref{lem5.1-2}, we derive
\begin{equation}
\|A^{\frac{\alpha}{2}} g(t_{0})\|^{2}+\varepsilon(t_{0})\|A^{\frac{1+\alpha}{2}} g(t_{0})\|^{2} \leqslant R_{\xi}(t),
\label{5.36-2}
\end{equation}
where $R_{\xi} (t)=C \rho_{\xi}^{2 p+2}(t) .$

As a conclusion of (\ref{5.36-2}), (\ref{5.9-2}) follows directly. $\hfill$$\Box$

\begin{Remark}\label{re5.1-2}
In Lemma $\ref{lem5.2-2}$ we assume $0<\sigma<\min \left\{\eta, \frac{-\varepsilon^{\prime}(t)}{\varepsilon(t)}, \frac{1}{\left(1+\lambda_{1}\right) \lambda_{1}^{-1}}\right\}$ with $0<\eta<m\lambda_{1}$ and in this lemma we assume that $0<\delta<\min \left\{(p+1) \sigma, \frac{16 m-2-\xi}{8\left(\lambda_{1}^{-1}+\varepsilon(t)\right)}\right\}$, then it follows the more accurate range of $\sigma :$
\begin{equation}
0<\sigma<\min \left\{m \lambda_{1}, \frac{-\varepsilon^{\prime}(t)}{\varepsilon(t)}, \frac{1}{\left(1+\lambda_{1}\right) \lambda_{1}^{-1}}, \frac{\delta}{p+1}\right\}.
\label{5.55-2}
\end{equation}
\end{Remark}
\begin{Lemma}\label{lem5.4-2}
For any $t_{0} \in \mathbb{R}$ and $\tau>0$, if $u_{0} $ is a subset of bounded set $B \subset \mathcal{H}_{t}(\Omega)$, the solution $u^{\xi}\left(t_{0}, t_{0}-\tau\right) u_{0}$ of problem $(\ref{1.1-2})$ converges to the solution $u(\tau) u_{0}$ of the following unperturbed equation of $(\ref{1.1-2}):$
\begin{equation}
\left\{\begin{array}{ll}
u_{t}-\varepsilon(t) \Delta u_{t}-a(l(u)) \Delta u=f(u) & {\rm{ in }} \,\, \Omega \times(\tau, \infty), \\
u=0 & {\rm{on}} \,\, \partial \Omega\times(\tau, \infty), \\
u(x, \tau)=u_{\tau}  & \, x \in \Omega,
\end{array}\right.\label{5.56-2}
\end{equation}
that is, \begin{equation}
\lim _{\xi \rightarrow 0^{+}} \sup _{u_{0} \in B}\|u^{\xi}(t_{0}, t_{0}-\tau) u_{0}-u(\tau) u_{0}\|_{\mathcal{H}_{t}}=0.
\end{equation}
\end{Lemma}
$\mathbf{Proof.}$ Let $\omega^{\xi}=u^{\xi}\left(t, t_{0}-\tau\right) u_{0}-u\left(t-t_{0}+\tau\right) u_{0}$, then $\omega^{\xi}$ satisfies
\[\left\{ {\begin{array}{*{20}{l}}
{w_t^\xi  - \varepsilon (t)\Delta w_t^\xi  - a(l({u^\xi }))\Delta {u^\xi } + a(l(u))\Delta u = f\left( {{u^\xi }} \right) - f(u) + \xi h(x,t),}&{{\rm{ in }}\,\,\Omega  \times (\tau ,\infty ),}\\
{\omega^\xi  = 0}&{{\rm{ on }}\,\,\partial \Omega  \times (\tau ,\infty ),}\\
{{\omega ^\xi }\left( {{t_0} - \tau } \right) = 0}&{{\mkern 1mu} {\mkern 1mu} x \in \Omega .}
\end{array}} \right.\]

Choosing $\omega^{\xi}$ as the test function of the above function in $L^{2}(\Omega)$, we obtain
\begin{equation}
\begin{aligned}
& \frac{d}{d t}(\|w^{\xi}\|^{2}+\varepsilon(t)\|\nabla w^{\xi}\|^{2})-\varepsilon^{\prime}(t)\|w^{\xi}\|^{2}+2 a(l(u^{\xi}))\|\nabla w^{\xi}\|^{2} \\
=& 2(a(l(u))-a(l(u^{\xi})))(\nabla u, \nabla w^{\xi})+2(f(u^{\xi})-f(u), w^{\xi})+\xi(h(x, t), w^{\xi}).
\end{aligned}
\label{5.57-2}
\end{equation}

Then the following inequalities follows by similar calculations to (\ref{3.25-2}) and (\ref{3.26-2}),
\begin{equation}
2(a(l(u))-a(l(u^{\xi})))(\nabla u, \nabla w^{\xi}) \leqslant 2 m\|\nabla u^{\xi}-\nabla u\|^{2}+\frac{(L a(R))^{2}\|l\|^{2}\|\nabla u\|^{2}\left\|u^{\xi}-u\right\|^{2}}{2 m},
\label{5.58-2}
\end{equation}
and
\begin{equation}
\begin{aligned}
2(f(u^{\xi})-f(u), \nabla w^{\xi}) & \leq C \int_{\Omega}(1+|u^{\xi}|^{p}+|u|^{p})|\omega^{\xi}|^{2} d x \\
& \leq C(\int_{\Omega}(1+|u^{\xi}|^{p}+|u|^{p})^{\frac{2 N}{(N-2) p}} d x)^{\frac{(N-2) p}{2 N}}(\int_{\Omega}|w^{\xi}|^{\frac{2 N}{N-2}} d x)^{\frac{N-2}{2 N}} \\
& \times(\int_{\Omega}|w^{\xi}|^{\frac{2 N}{2N-(N-2)(p+1)}} d x)^{\frac{2N-(N-2)(p+1)}{2 N}} \\
& \leq C(1+\|A^{\frac{1}{2}} u^{\xi}\|^{p}+\|A^{\frac{1}{2}} u\|^{p})\|A^{\frac{1}{2}} w^{\xi}\|\|A^{\frac{1-\alpha}{2}} w^{\xi}\| \\
&\leq C(1+\|A^{\frac{1}{2}} u^{\xi}\|^{p}+\|A^{\frac{1}{2}} u\|^{p})(\|w^{\xi}\|^{2}+\varepsilon(t)\|\nabla w^{\xi}\|^{2}) \\
& \leq C (1+ \|A^{\frac{1}{2}} u^{\xi} \|^{\frac{4}{N-2}}+ \|A^{\frac{1}{2}} u \|^{\frac{4}{N-2}} )\|w^{\xi}\|_{\mathcal{H}_{t}}^{2}.
\end{aligned}
\label{5.59-2}
\end{equation}

By the Young inequality, we can derive
\begin{equation}
2(\xi h(x, t), \omega^{\xi}) \leqslant \frac{\xi^{2}}{\lambda_{1}}\|h(x, t)\|^{2}+\lambda_{1}\|\omega^{\xi}\|^{2}.
\label{5.60-2}
\end{equation}

Inserting $(\ref{5.58-2})-(\ref{5.60-2})$ into (\ref{5.57-2}), from (\ref{1.4-2}) and $\varepsilon(t)$ is a decreasing function, then using the Sobolev embedding theorem, we conclude
\begin{equation}
\begin{aligned}
\frac{d}{d t}(\|w^{\xi}\|_{\mathcal H_{t}}^{2}) &=\frac{d}{d t}(\|w^{\xi}\|^{2}+\varepsilon(t)\|\nabla w^{\xi}\|^{2})\\
&\leqslant C(1+\|A^{\frac{1}{2}} u^{\xi}\|^{\frac{4}{N-2}}+\|A^{\frac{1}{2}} u\|^{\frac{4}{N-2}})\|w^{\xi}\|^{2}_{\mathcal H_{t}}+C {\xi}^{2}\|h(x, t)\|^{2}.
\end{aligned}
\label{5.61-2}
\end{equation}

By the Gronwall inequality, we obtain
\begin{equation}
\|w^{\xi}(t_{0})\|_{\mathcal H_{t}}^{2} \leqslant C \xi^{2} e^{\int_{t_{0}-\tau}^{t_{0}}(1+\|A^{\frac{1}{2}} u^{\xi}\|^{\frac{4}{N-2}}+\|A^{\frac{1}{2}} u\|^{\frac{4}{N-2}}) d t} \cdot \int_{t_{0}-\tau}^{t_{0}}\left\|h\left(x,t\right)\right\|^{2} d t.
\label{5.62-2}
\end{equation}

Then the desired result follows directly. $\hfill$$\Box$

The main result of this section is as follows.
\begin{Theorem}\label{th5.1-2}
Under the assumptions of Theorems $\ref{th3.1-2}-\ref{th3.2-2}$ and Lemmas $\ref{lem5.1-2}-\ref{lem5.4-2}$, assume further that the function $h(x, t)$ satisfies $(\ref{4.11-2})$, then the pullback attractor $\mathcal{A}_{\xi}=\{{A}_{\xi}(t)\}_{t \in \mathbb R}$ of problem $(\ref{1.1-2})$ with $\xi >0$ and the global attractor $A$ for problem $(\ref{1.1-2})$ with $\xi =0$ satisfy
\begin{equation}
\lim _{\xi \rightarrow 0^{+}} \operatorname{dist}_{\mathcal{H}_{t}}\left(A_{\xi}(t), A\right)=0,
\label{5.01-2}
\end{equation}
for all $t \in \mathbb R$.
\end{Theorem}
$\mathbf{Proof.}$ From Lemma $\ref{lem5.1-2}$ to Lemma $\ref{lem5.4-2}$, we deduced that Lemma \ref{1m2.2-2} holds. Then by Lemma \ref{th2-2}, (\ref{5.01-2}) follows directly. $\hfill$$\Box$

$\mathbf{Acknowledgment}$

This paper was partially supported by the NNSF of China with contract number 12171082 and by the fundamental funds for the central universities with contract number $2232022G$-$13$.

$\mathbf{Conflict\,\,of\,\,interest\,\,statement}$

The authors have no conflict of interest.


\newpage
\begin{thebibliography}{lllp}
\setlength{\itemsep}{- 2mm}
\bibitem{acr.2} M. Anguiano, T. Caraballo and J. Real, An exponential growth condition in $H^{2}$ for the pullback attractor of a non-autonomous reaction-diffusion equation, Nonlinear Anal., 72(11)(2010), 4071-4075.
\bibitem{ab1.2} C. T. Anh and T. Q. Bao, Pullback attractors for a class of non-autonomous nonclassical diffusion equations, Nonlinear Anal., 73(2)(2010), 399-412.
\bibitem{ab2.2} C. T. Anh and T. Q. Bao, Dynamics of non-autonomous nonclassical diffusion equations on $\mathbb R^{n}$, Comm. Pure Applied Anal., 11(3)(2012), 1231-1252.
\bibitem{bv.2} A. V. Babin, M. I. Vishik, Attractors of evolution equations,  Elsevier Science Publishers B.V., Netherlands, 1992.
\bibitem{cc} T. Caraballo and D. N. Cheban, On the structure of the global attractor for non-autonomous difference equations with weak convergence, J. Diff. Equ. Appl., 18(4)(2012), 535-551.
\bibitem{chm2} T. Caraballo, M. Herrera-Cobos and P. Mar\'{i}n-Rubio, Long-time behavior of a non-autonomous parabolic equation with nonlocal diffusion and sublinear terms, Nonlinear Anal., 121(7)(2015), 3-18.
\bibitem{chm16} T. Caraballo, M. Herrera-Cobos and P. Mar\'{i}n-Rubio, Robustness of nonautonomous attractors for a family of nonlocal reaction-diffusion equations without uniqueness, Nonlinear Dyn., 81(1)(2016), 35-50.
\bibitem{chm3.2} T. Caraballo, M. Herrera-Cobos and P. Mar\'{i}n-Rubio, Asymptotic behaviour of nonlocal $p$-Laplacian reaction-diffusion problems, J. Math. Anal. Appl., 459(2)(2018), 997-1015.
\bibitem{chm4.2} T. Caraballo, M. Herrera-Cobos and P. Mar\'{i}n-Rubio, Robustness of time-dependent attractors in $H^{1}$-norm for nonlocal problems, Disc. Contin. Dynam. Syst., 23(3)(2018), 1011-1036.
\bibitem{chm} T. Caraballo, M. Herrera-Cobos and P. Mar\'{i}n-Rubio, Time-dependent attractors for nonautonomous nonlocal reaction-diffusion equations, Proc. Roy. Soc. Edinburgh-A, 148(5)(2018), 957-981.
\bibitem{clr.2} T. Caraballo, J. A. Langa and J. C. Robinson, Upper semicontinuity of attractors for small random perturbations of dynamical systems, Comm. Partial Differ. Equ., 23(9-10)(1998), 1557-1581.
\bibitem{cjr1.2} A. N. Carvalho, J. A. Jos$\acute{e}$ and J. C. Robinson, On the continuity of pullback attractors for evolution processes, Nonlinear Anal., 71(5-6)(2009), 1812-1824.
\bibitem{cv} V. V. Chepyzhov and M. I. Vishik, Attractors for equations of mathematical physics, American Mathematical Society, America, 2002.
\bibitem{ct2.2} C. A. Cung and N. D. Toan, Existence and upper semicontinuity of uniform attractors in $H^{1}(\mathbb R^{N}) $ for nonautonomous nonclassical diffusion equations, Annales Polonici Math., 111(3)(2014), 271-295.
\bibitem{ct1.2} C. A. Cung and N. D. Toan, Uniform attractors for non-autonomous nonclassical diffusion equations on $\mathbb R^{N}$, Bulletin of the Korean Math. Soc., 51(5)(2014), 1299-1324.
\bibitem{Evans} L. C. Evans, Partial  differential  equations, American Mathematical Society, America, 1998.
\bibitem{gmr.2} J. Garc\'{i}a-Luengo, P. Mar\'{i}n-Rubio and J. Real, Pullback attractors in $V$ for non-autonomous 2$D$-Navier-Stokes equations and their tempered behaviour, J. Diff. Equ., 252(8)(2012), 4333-4356.
\bibitem{gw.2} B. Guo and B. Wang, Upper semicontinuity of attractors for the reaction diffusion equation, Acta Math. Sci., 18(2)(1998), 139-145.
\bibitem{Lions} J. L. Lions, Quelques m$\acute{e}$thodes de r$\acute{e}$solutions des problèms aus limites nonlin$\acute{e}$aries, Dunod Gauthier-Villars, France, 1969.
\bibitem{lm} J. L. Lions and E. Magenes, Non-homogeneous boundary value problem and applications, Springer-Verlag, Germany, 1972.
\bibitem{mwx} Q. Ma, X. Wang and L. Xu, Existence and regularity of time-dependent global attractors for the nonclassical reaction-diffusion equations with lower forcing term, Boundary Value Problems, 10(1)(2016), 1-11.
\bibitem{ps.2} V. Pata and M. Squassina, On the strongly damped wave equation, Comm. Math. Phys. 253(3)(2005), 511-533.
\bibitem{psz} X. Peng, Y. Shang and X. Zheng, Pullback attractors of nonautonomous nonclassical diffusion equations with nonlocal diffusion, Zeitschrift f\"ur Angewandte Mathematik und Physik, 69(4)(2018), 1-14.
\bibitem{Q1} Y. Qin, Integral and discrete inequalities and their applications, Vol I: Springer International Publishing AG, Switzerland, 2016.
\bibitem{Q2} Y. Qin, Integral and discrete inequalities and their applications, Vol II: Springer International Publishing AG, Switzerland, 2016.
\bibitem{Q3} Y. Qin, Analytic inequalities and their applications in PDEs, Birkhauser Verlag AG, Switzerland, 2017.
\bibitem{qy.2} Y. Qin and B. Yang, Existence and regularity of time-dependent pullback attractors for the non-autonomous nonclassical diffusion equations, Proc. Roy. Soc. Edinburgh-A, 2021, DOI: https://doi.org/10.1017/prm.2021.65.
\bibitem{r} J. C. Robinson, Infinite-dimensional dynamical systems, Cambridge University Press, England, 2011.
\bibitem{sy} C. Sun and M. Yang, Dynamics of the nonclassical diffusion equations, Asymptotic Anal., 59(1-2)(2008), 51-81.
\bibitem{w.2} Y. Wang, On the upper semicontinuity of pullback attractors for multi-valued processes, Quarterly Appl. Math., 71(2)(2013), 369-399.
\bibitem{wlq.2} Y. Wang, P. Li and Y. Qin, Upper semicontinuity of uniform attractors for nonclassical diffusion equations, Boundary Value Problems, 84(7)(2017), 1-11.
\bibitem{wq.2} Y. Wang and Y. Qin, Upper semicontinuity of pullback attractors for nonclassical diffusion equations, J. Math. Phy., 51(2)(2010), 1-12.
\bibitem{ww.2} Y. Wang and L. Wang, Trajectory attractors for nonclassical diffusion equations with fading memory, Acta Math. Sci., 33(3)(2013), 721-737.
\bibitem{wwq.2} L. Wang, Y. Wang and Y. Qin, Upper semicontinuity of attractors for nonclassical diffusion equations in $H^{1}(\mathbb R^{3})$, Appl. Math. Comput., 240(1)(2014), 51-61.
\bibitem{xz.2} Y. Xie and K. Zhu, Upper semicontinuity of uniform attractors for nonclassical diffusion equations, Advances in Diff. Equ., 75(1)(2021), 1-17.
\bibitem{yqlm.2} X. Yang, Y. Qin, Y. Lu and T. Ma, Dynamics of 2D incompressible non-autonomous Navier-Stokes equations on Lipschitz-like domains, Appl. Math. Opt., 83(3)(2021), 2129-2183.
\bibitem{zb.2} F. Zhang and L. Bai, Pullback attractors in the weighted space for multi-valued process generated by the non-autonomous nonclassical diffusion equations with unbounded delays without uniqueness of solutions, Appl. Analy., 99(8)(2020), 1436-1452.
\bibitem{zl.2} F. Zhang and Y. Liu, Pullback attractors in $H^{1} (\mathbb R^{N})$ for non-autonomous nonclassical diffusion equations, Dyna. Syst., 29(1)(2014), 106-118.
\bibitem{zs2} K. Zhu and C. Sun, Pullback attractors for nonclassical diffusion equations with delays, J. Math. Phy., 56(9)(2015), 1-20.
\bibitem{zxz} K. Zhu, Y. Xie and F. Zhou, Attractors for the nonclassical reaction-diffusion equations on time-dependent spaces, Boundary Value Problems, 95(1)(2020), 1-14.
\end{thebibliography}
\end{document}